\documentclass[letterpaper]{ijmart}
\usepackage{hyperref}
\usepackage[utf8]{inputenc}
\usepackage{amsmath}
\usepackage{amsfonts}
\usepackage{amssymb}
\usepackage{diagrams}
\usepackage{stmaryrd}
\usepackage{eucal}
\usepackage{graphicx}
\usepackage{indentfirst}

\usepackage{mathrsfs}

\usepackage[T1]{fontenc}

\normalfont

\newtheorem{thm}{Theorem}[section]
\newtheorem{prop}[thm]{Proposition}
\newtheorem{lemma}[thm]{Lemma}
\newtheorem{corr}[thm]{Corollary}

\newtheorem{introthm}{Theorem}%
\newtheorem{introprop}[introthm]{Proposition}%
\newtheorem*{prop*}{Proposition}

\theoremstyle{definition}
\newtheorem{opr}[thm]{Definition}

\theoremstyle{remark}
\newtheorem{rem}[thm]{Remark}
\newtheorem{ntn}[thm]{Notation}

\makeatletter
\@addtoreset{equation}{section}
\makeatother

\newenvironment{introequation}{%

\begin{equation}%
}%
{\end{equation}%
	}%

\newcommand{\bC}{{\mathbb C}}
\newcommand{\bR}{{\mathbb R}}

\newcommand{\bD}{{\mathbb D}}

\newcommand{\bI}{{\mathbb I}}

\newcommand{\bL}{{\mathbb L}}

\newcommand{\cB}{{\mathcal B}}
\newcommand{\cA}{{\mathcal A}}
\newcommand{\cC}{{\mathcal C}}
\newcommand{\cD}{{\mathcal D}}

\newcommand{\cF}{{\mathcal F}}

\newcommand{\cN}{{\mathcal N}}
\newcommand{\cM}{{\mathcal M}}
\newcommand{\cL}{{\mathcal L}}
\newcommand{\cZ}{{\mathcal Z}}

\newcommand{\cS}{{\mathcal S}}

\newcommand{\cW}{{\mathcal W}}

\newcommand{\cX}{{\mathcal X}}
\newcommand{\cY}{{\mathcal Y}}
\newcommand{\cI}{{\mathcal I}}
\newcommand{\cE}{{\mathcal E}}

\newcommand{\cR}{{\mathcal R}}

\newcommand{\cO}{{\mathcal O}}

\newcommand{\h}{{\mathsf h}}

\renewcommand{\phi}{\varphi}

\newarrow {Into} {boldhook}--->
\newarrow {Mapsto} |--->
\newcommand{\Set}{{\rm \bf Set}}
\newcommand{\Fin}{\Gamma}

\newcommand{\SSet}{{\rm \bf SSet}}
\newcommand{\Cat}{{\rm \bf Cat }}

\newcommand{\Top}{{\rm \bf Top}}

\newcommand{\DVect}{{\rm \bf DVect}}
\newcommand{\DGMod}{{\rm \bf DVect}}

\newcommand{\Ob}{{\rm Ob} \,}

\newcommand{\Ho}{{\rm Ho} \,}

\newcommand{\bT}{\mathbb T}

\newcommand{\bfE}{\mathbf E}
\newcommand{\bF}{\mathbb F}

\newcommand{\op}{{\sf op}}
\newcommand{\bd}{\mathbf d}
\newcommand{\bc}{\mathbf c}

\renewcommand{\lim}{\operatorname{\varprojlim}}
\newcommand{\colim}{\operatorname{\varinjlim}}
\newcommand{\hocolim}{\operatorname{hocolim}}

\newcommand{\Sect}{{\sf Sect}}
\newcommand{\PSect}{{\sf PSect}}
\newcommand{\DSect}{{\sf DSect}}

\newcommand{\Cart}{{\sf Cart}}

\newcommand{\Mat}{\mathscr{Mat}}
\newcommand{\sR}{\mathscr{R}}
\newcommand{\sL}{\mathscr{L}}
\newcommand{\sA}{\mathscr{A}}
\newcommand{\sS}{\mathscr{S}}
\newcommand{\sG}{\mathscr{G}}
\newcommand{\sD}{\mathscr{D}}

\newcommand{\Kappa}{\mathrm{K}}

\newcommand{\rE}{{\mathrm E}}
\newcommand{\rF}{{\mathrm F}}

\newcommand{\rA}{{\mathrm A}}
\newcommand{\rO}{{\mathrm O}}

\newcommand{\Fun}{\operatorname{\mathsf{Fun}}}

\author{Edouard Balzin}
\title{The formalism of Segal sections}
\date{}

\makeatletter
\ifcsname phantomsection\endcsname
    \newcommand*{\qrr@gobblenexttocentry}[5]{}
\else
    \newcommand*{\qrr@gobblenexttocentry}[4]{}
\fi
\newcommand*{\addsubsection}{%
    \addtocontents{toc}{\protect\qrr@gobblenexttocentry}%
    \subsection}
\makeatother

\begin{document}
\begin{abstract}
Given a family of model categories $\cE \to \cC$, we associate to it a
homotopical category of derived, or Segal, sections $\DSect(\cC,\cE)$ that models the higher-categorical
sections of the localisation $L\cE \to \cC$. The derived sections provide an
alternative, strict model for various higher algebra objects appearing in the work of Lurie.
We prove a few results concerning the properties of the homotopical category $\DSect(\cC,\cE)$,
and as an example, study its behaviour with respect to the base-change along a select class of functors.
\end{abstract}
	\maketitle
    \thispagestyle{empty}

\tableofcontents

\section*{Introduction}
\addtocontents{toc}{\protect\setcounter{tocdepth}{1}}
\subsection*{Segal objects} The formalism presented in this paper was developed in the study of homotopy algebraic
structures as described by Segal and generalised by Lurie. We begin the introduction
by describing this context.

Denote by $\Fin$ the category whose objects are finite sets and morphisms
are given by partially defined set maps. Each such morphism between $S$ and $T$
can be depicted as $S \supset S_0 \to T$. A \emph{$\Fin$-space} is simply a functor
$X: \Fin \to \Top$ taking values in the category of topological spaces. Each such functor
can be evaluated on the maps in $\Fin$ that have the form $S \supset {s} \to 1$, where
$s \in S$ and $1$ is the one-element set. As a result, for each $S \in \Fin$ we get the map
\begin{introequation}
\label{introsegmap}
X(S) \longrightarrow X(1)^S
\end{introequation}%
and the $\Fin$-space is called \emph{Segal} \cite{SEG} if the map (\ref{introsegmap}) is a homotopy
equivalence for all $S \in \Fin$.

For each $S$, we also have the map $S \stackrel = \supset S \to 1$, and evaluating $X$
on it, combined with (\ref{introsegmap}), gives the following span in $\Top$:
\begin{introequation}
\label{introsegdia}
X(1)^S \stackrel \sim \longleftarrow X(S) \longrightarrow X(1);
\end{introequation}%
choosing an inverse equivalence to the left map, we get a non-canonical operation $X(1)^S \to X(1)$.
There is no ambiguity in choice of the inverse in the homotopy category of topological spaces
$\Ho \Top$: the result is a commutative monoid structure on $X(1)$ in the homotopy category.
The full structure that $X(1)$ carries is that of an $\rE_\infty$, or a homotopy commutative monoid.

Instead of $\Fin$ one can work with other categories. For example, take the category
$\Delta$ which consists of finite categories $[n]= 0 \to 1 \to ... \to n$ ($n \geq 0$)
and functors between them.
We can then consider functors $Y:\Delta^\op \to \Top$ such that
$Y([0]) = *$ and for each $n\geq 1$, the map $Y([n]) \to Y([1])^n$, induced by all the
consecutive interval inclusions $[1] \to [n]$, is a homotopy equivalence. Such Segal
$\Delta$-spaces are known to describe homotopy associative, or $\rE_1$-monoids. The
functor $\Delta^\op \to \Fin$ which maps $[n]$ to the set of its $n$ generating arrows,
allows us to pull $\Gamma$-spaces back to $\Delta^\op$: every commutative monoid is also
an associative one.
There are further examples of small categories that parametrise operations up to homotopy.
For instance, \cite{BAR} introduces the notion of an operator category and gives
examples of categories that model $\rE_m$-algebras. The case $m=2$ can also be approached
geometrically using the exit path category \cite{WOOLF} of the $2$-disk. One can also
introduce categories that permit to describe an $\rE_m$-algebra over another algebra. All such categories
$\rO$ come with a functor $\rO \to \Fin$, so that we can pull Segal $\Gamma$-spaces back to $\rO$,
reflecting the ``initial'' property of homotopy commutative structures.

Let us attempt to replace $\Top$ with another homotopical category (by which we simply mean
a category with a distinguished class of weak equivalences) possessing a monoidal
structure. For example, take the category $\DVect_k$ of chain complexes of vector spaces over a field $k$.
The maps like (\ref{introsegmap}) will take the form $X(S) \to \oplus_S X(1)$, which
is not very interesting for the purposes of algebra. Indeed, the tensor product $\otimes_k$
is not Cartesian. The observation of \cite{LU} addresses this issue by making the following construction.

\subsection*{Algebras as sections} Given any symmetric monoidal category $(\cM, \otimes, \bI)$, let us define the category $\cM^\otimes$. Its
objects are pairs $(S, \{X_s\}_{s \in S})$ where $S \in \Fin$ and each $X_s$ is an object of $\cM$.
A morphism $(S, \{X_s\}_{s \in S}) \to (T, \{Y_t\}_{t \in T})$ consists of a partially
defined map $f:S \to T$, and for each $t \in T$, of a morphism $\otimes_{s \in f^{-1}(t)} X_s \to Y_t$.
When $f^{-1}(t)$ is empty, the monoidal product over it equals $\bI$.
The compositions can then be defined with the help of the coherence isomorphisms
for the product $\otimes$ and the unit object.

The forgetful functor $p: \cM^\otimes \to \Fin$ is a \emph{Grothendieck opfibration} \cite{SGA1}:
the assignment $S \mapsto \cM^\otimes(S):= p^{-1}(S) = \cM^S $ is functorial, but in a weak way. Each map $f: S \to T$
corresponds to a functor $f_!:\cM^\otimes(S) \to \cM^\otimes(T)$, sending $(S, \{X_s\}_{s \in S})$ to
$(T, \{(\otimes_{s \in f^{-1}(t)} X_s)\}_{t \in T})$. We see that for composable maps
$f,g$, one has $g_! f_! \cong (gf)_!$ and such isomorphisms, while suitably coherent, are nonetheless
not equalities. It is more convenient to work with Grothendieck (op)fibrations since
this allows to keep all these isomorphisms implicit.

Denoting for simplicity two objects of $\cM^\otimes$ as $X,Y$, any map $\alpha: X \to Y$
factors as $X \to f_! X \to Y$ where $f = p (\alpha)$; the map $\alpha$ is called \emph{cocartesian}
if the induced map $f_! X \to Y$ is an isomorphism. Consider now a section of $p$, that is
a functor $\cA: \Fin \to \cM^\otimes$ with $p \cA = id$. Demand further that given any \emph{inert} map
in $\Fin$, that is a map of the form $S \supset T \stackrel = \to T$, the induced
map $\cA(S) \to \cA(T)$ is cocartesian. It is elementary to verify that each $\cA(S)$
is isomorphic to $S$ copies of $\cA(1) \in \cM^\otimes(1) = \cM$. The remaining information provided by the section $\cA$ equips $\cA(1)$
with the structure of a commutative monoid in $\cM$: the value of $\cA$ on the map
$S \stackrel = \supset S \to 1$ gives multiplication $\cA(1)^S \to \cA(1)$, the composition
property insures that everything is determined by this map for $|S|=2$ (which is forced to be associative)
and the action of automorphisms in $\Fin$ forces the multiplication to be commutative.

Setting $\cM = \DVect_k$, it is well known that the homotopy theory of commutative
$dg$-algebras is an inadequate model for $\rE_\infty$-algebras in $\DVect_k$ if
$char k \neq 0$.  From the fibrational perspective, the reason for this is that
the transition functors of the opfibration $\cM^\otimes \to \Fin$ lack adjoints
(they are multifold tensor products). One cannot thus apply the model structure results of
\cite{BARLEFT} or pull back
to a suitable Reedy category and use \cite{BALRMSF}.
Another problem arises if we try different parametrising categories $\rO$ described
above. We can consider functors $\rO \to \cM^\otimes$ that commute with the projection to
$\Fin$, subject to some cocartesian conditions. This way one can indeed obtain associative algebra
in the case of $\Delta^\op \to \Fin$, but the operator categories parametrising $\rE_m$-structures
will give commutative algebras for $m \geq 2$, due to the Eckmann-Hilton argument.

In order to deal with the problems described above, \cite{LU} passes to
the world of symmetric monoidal infinity-categories,
and performs all the constructions on that level. More precisely the machinery of \cite{LU}
is developed for the infinity-operads,
and is not immediately applicable to arbitrary operation-indexing categories $\rO \to \Fin$ some of which \cite{BALTH} have no operadic
analogues (technically speaking, one can attempt taking a fibrant replacement of $\rO$ in the category of infinity-operads,
but this will result in the loss of control over the combinatorics of $\rO$).
This approach is also, arguably, quite
technical: the chapters of \cite{LU} devoted to colimits of algebras and adjoint functors
contain heavy combinatorics necessary even for colimits of algebras in the presentable setting (in the sense of \cite[5.5]{LUHTT}).
These issues motivated us to look for an alternative approach.

\subsection*{Segal sections} Let us pass from $\cM^\otimes \to \Gamma$ to more general
Grothendieck opfibrations $\cE \to \cC$ over small categories $\cC$. For each $c \in \cC$
we thus have the fibre category $\cE(c)$, and each map $f: c \to c'$ produces a transition
functor $f_!:\cE(c) \to \cE(c')$ (the compositions are respected up to an isomorphism).
Define $\Sect(\cC,\cE)$
to be the category of sections of $\cE \to \cC$; it is defined as the fibre of the post-composition
functor $\Fun(\cC,\cE) \to \Fun(\cC,\cC)$ (with $\Fun$ denoting the functor categories) over the identity functor.
Any section $X: \cC \to \cE$ associates to a map $f: c \to c'$ the map $f_! X(c) \to X(c')$ in $\cE(c')$;
if the latter map is an isomorphism, the section $X$ is cocartesian along $f$. We can
impose the cocartesian condition along a subset $\cS$ of maps of $\cC$, and denote
the resulting subcategory as $\Sect_\cS(\cC,\cE)$. This category is a generalisation
of the algebra objects category in the case of $\cM^\otimes \to \Fin$.

Assume now that each $\cE(c)$ has weak equivalences $\cW(c)$.
For simplicity we ask here that the transition functors $f_!$
preserve weak equivalences. We can then take the higher-categorical localisation
of $\cE$ with respect to all $\cW(c)$, and as observed in \cite{HINICH,MGEE1},
the resulting infinity-functor $L \cE \to \cC$ is a cocartesian fibration in infinity-categories
as defined in \cite{LUHTT}. One can use the language of \cite{LUHTT} to define the higher categories
$\Sect(\cC,L\cE)$
and $\Sect_\cS(\cC,L\cE)$. We would like to understand the properties of these higher categories
as a function of the properties of $\cE$, and ideally have convenient strict models for them.

Our answer to this problem is the notion of a \emph{derived}, or \emph{Segal}, section (we will use
both names interchangeably
in this text).
Informally, such an object $X$
associates to $c \in \cC$ an object $X(c) \in \cE(c)$, and to each map $f: c \to c'$,
a diagram
\begin{introequation}
\label{introsegalsectdia}
f_! X(c) \stackrel{ \sim }{\longleftarrow} X_f \longrightarrow X(c')
\end{introequation}%
with $f_!: \cE(c) \to \cE(c')$ being the transition functor and the left map in
(\ref{introsegalsectdia}) belonging to $\cW(c)$. If we were to invert it, we would
get a map $f_! X(c) \to X(c')$: in this sense, a Segal section is a derived version
of the concept of a section.
There are also additional data that ensure the coherence
with respect to the compositions in $\cC$. Finally, if $f \in \cS$, then both
maps in the diagram (\ref{introsegalsectdia}) are required to be weak equivalences.

Below we explain more formally our definition of Segal sections;
those readers that are only interested in the properties of the homotopical category of derived
sections are invited to skip further.

The start is to replace the base category $\cC$. The well-known construction associates to
it the category $\Delta / \cC$, whose objects are functors $[n] \to \cC$ (in other words,
strings of arrows $c_0 \to ... \to c_n$), and the morphisms are functors $[n] \to [m]$ over $\cC$.
The \emph{simplicial replacement} $\bC$ of $\cC$ is defining by simply putting
$\bC := (\Delta / \cC)^\op$.

Among the maps in $\Delta$, we can consider those interval inclusions $[m] \hookrightarrow [n]$
that send $0$ to $0$. The maps in $\bC$ with such underlying $\Delta$-map look like
$$
(c_0 \to ... \to c_n) \longrightarrow (c_0 \to ... \to c_m)
$$
and are called Segal maps in this paper. The reason for this terminology is the known
fact that the functor $h:\bC \to \cC$ sending $c_0 \to ... \to c_n$ to $c_0$, is the
(higher-categorical) localisation of $\bC$ along Segal maps (Proposition \ref{locprop}).
As a result, given a category $\cM$ with weak equivalences $\cW$, we can represent
the functors from $\cC$ to the localisation of $\cM$ with respect to $\cW$ (in ordinary or higher
sense) as functors from $\bC$ to $\cM$ sending the Segal maps to $\cW$. Given
such functor $X$, its evaluation on the span $c_0 \leftarrow (c_0 \stackrel f \to c_1) \to c_1$
gives a diagram
\begin{introequation}
X(c_0) \stackrel{ \in \cW}{\longleftarrow} X(c_0 \stackrel f \to c_1) \longrightarrow X(c_1)
\end{introequation}%
that can be compared with (\ref{introsegdia}) and is a version of (\ref{introsegalsectdia})
when no non-trivial opfibration is involved. In this way, $X: \bC \to \cM$ can be viewed as
a ``Segal functor'' from $\cC$ to $\cM$.

The next step of the construction concerns the opfibration $\cE \to \cC$. Its extension
to $\bC$ can be done using the functor $h:\bC \to \cC$. While this technique is useful when dealing with
Quillen presheaves \cite{BALRMSF}, the result is not illuminating
in the case of $\cM^\otimes \to \Gamma$ due to the already mentioned lack of adjoints.
Instead, one has to turn to the \emph{transpose}
fibration $\cE^\top \to \cC^\op$ (Definition \ref{transposefib}): it has the same fibres $\cE(c)$, and has the same
transition functors inherited from $\cE \to \cC$, but read in a contravariant direction
over $\cC^\op$. Now, any object $c_0 \to ... \to c_n \in \bC$, viewed as a functor $\bc:[n] \to \cC$,
allows us to write
$$
\bfE(\bc):= \Sect([n], \bc^{\op,*}\cE^\top),
$$
in other words, we consider the sections of the transpose fibration over the string
$c_0 \to ... \to c_n$. Moreover the assignment $\bc \mapsto \bfE(\bc)$ defines a functor: for each map $\alpha:\bc \to \bc'$
in $\bC$, we have the restriction $\alpha_!:\bfE(\bc) \to \bfE(\bc')$. Using the Grothendieck construction
\cite{VIST} we can take the associated opfibration denoted $\bfE \to \bC$.

Given an object $\bc=c_0 \to ... \to c_n$, the functor $\bfE(\bc) \to \bfE(c_0)=\cE(c_0)$
of evaluation at $c_0$ admits a right adjoint that induces an equivalence
$\bfE(c_0) \cong \Cart([n], \bc^{\op,*}\cE^\top) \subset \bfE(\bc)$ of $\bfE(c_0)$ with the full subcategory
of $\Sect([n], \bc^{\op,*}\cE^\top)$ consisting of sections taking all maps of $[n]$ to cartesian maps. A similar
situation happens along any Segal map $\alpha:\bc \to \bc'$,
and the existing adjunctions $\alpha_!: \bfE(\bc) \rightleftarrows \bfE(\bc'): \alpha^*$
mean that $\bfE \to \bC$ is also a fibration along Segal maps.

The following is Proposition {\ref{sectinpsect}} in the main text.

\begin{introprop}
\label{introcomparisonbasic}
The category $\Sect(\cC,\cE)$ is equivalent
to the full subcategory of $\Sect(\bC,\bfE)$ consisting of those $X$ such that for each Segal
map $\alpha:\bc \to \bc'$, the image $X(\alpha)$ is cartesian, meaning $X(\bc) \stackrel \sim \to \alpha^* X(\bc')$.
\end{introprop}
If $\cE \to \cC$ has the weak equivalences as before, each category $\bfE(\bc)$ also comes with weak equivalences, defined value-by-value. The category of \emph{derived sections}
$\DSect(\cC,\cE)$ is then defined as the full subcategory of $\Sect(\bC,\bfE)$ consisting of all $X$
such that each Segal $\alpha:\bc \to \bc'$ induces a weak equivalence
 $X(\bc) \stackrel \cW \longrightarrow \alpha^* X(\bc')$ (more generally, in the precedent
 expression
we will replace $\alpha^*$ by its right derived functor $ \bR \alpha^*$). The weak equivalences
 of derived sections are defined value-by-value, giving a homotopical structure on $\DSect(\cC,\cE)$.
Without further explanation, we also mention that one can take the subset $\cS$ of maps of $\cC$
into account (Definition \ref{deflocconstdersect}), resulting in a full homotopical
subcategory $\DSect_\cS (\cC,\cE) \subset \DSect(\cC,\cE)$ of $\cS$\emph{-locally constant} derived sections.



\subsection*{Summary of properties} The summary of the discussion above is that given
an opfibration $\cE \to \cC$, there exists another opfibration $\bfE \to \bC$ that
inherits weak equivalences from $\cE$, and the homotopical category of derived sections $\DSect(\cC,\cE)$
is defined as a full subcategory of $\Sect(\bC,\bfE)$ using a weakened version of cartesian section condition.
Let us now describe some of its properties.

We first study the case when
the fibres of $\cE \to \cC$ are model categories, and the transition functors preserve
fibrations and trivial fibrations (but not necessarily have adjoints): such opfibrations
are called \emph{model} in this text. One important example is the case of $\DVect_k^\otimes \to \Fin$ with $k$ any field, whose
infinity-sections model $\rE_\infty$-algebras over $k$ (one can use operation-indexing categories to include
$\rE_n$-algebras over $k$ and other structures). Another example is
 $\DVect_k^{\op,\otimes} \to \Fin$, corresponding to the monoidal structure on the
 opposite category of chain complexes over $k$: studying infinity-sections
 in this case corresponds to higher coalgebra objects in $L \DVect_k$.

The following proposition is the summary of subsection \ref{homcatofsegsect}:
\begin{introprop}
\label{intropropositionbasic}
Let $\cE \to \cC$ be a model opfibration. Then the functor $\bfE \to \bC$ is a bifibration in
model categories and Quillen pairs, its sections $\Sect(\bC,\bfE)$ is a model category for the
Reedy structure of \cite{BALRMSF},
and the full subcategory $\DSect(\cC,\cE) \subset \Sect(\bC,\bfE)$ is closed under homotopy
limits.

Moreover when the model categories $\cE(c)$ are combinatorial and the transition functors $f_!$
are accessible, there is a Bousfield localisation of the Reedy model structure on $\Sect(\bC,\bfE)$
whose fibrant objects are satisfy the derived section condition. In particular, there exist homotopy
colimits of derived sections.
\end{introprop}

If the transition functors along the maps from the subset $\cS$ preserve limits,
then the inclusion $\DSect_\cS(\cC,\cE) \subset \Sect(\bC,\bfE)$ is also stable under limits,
yet the existence of a localised model structure is not guaranteed (the situation is similar
to that encountered in \cite{BARRIGHT}). Nonetheless Proposition \ref{intropropositionbasic}
allows to establish the following important result. Given a model opfibration $\cE \to \cC$,
denote by $\cE_f \subset \cE$ the subopfibration spanned by the fibrewise-fibrant objects.
We can apply the infinity-localisation of \cite{HINICH} to get a cocartesian fibration in quasicategories
$L \cE_f \to \cC$. The following result is a consequence of the comparison theorem
given in \cite{BALRMSF} and is a combination of Theorem \ref{thmcomparisonforlocconstdersect} and Corollary
\ref{corrlocdsectarepresentable}
in the main text:

\begin{introthm}
\label{introthmcomparisonmodel}
Let $\cE \to \cC$ be a model opfibration and $\cS$ a (possibly empty) subset of maps of $\cC$.
The infinity-category $\Sect_\cS(\cC,L \cE_f)$ consisting of sections $\cC \to L\cE_f$
that send the maps of $\cS$ to cocartesian maps of $L \cE_f$ is equivalent
to the infinity-localisation $L \DSect_\cS(\cC,\cE)$, in a manner compatible with the base-change.

Assume
that the fibres of $\cE \to \cC$ are combinatorial model categories, the transition functors are accessible,
and the transition functors along the maps in $\cS$ preserve limits. Then both
$L \DSect_\cS(\cC,\cE)$ and $\Sect_\cS(\cC,L \cE_f)$ are presentable infinity-categories.
\end{introthm}

We remark that even in the non-homotopical context, it is not obvious why the sections
of $\cE \to \cC$ have colimits, as they are indeed not given fibrewise. One can view the
idea of a Segal section to be a smarter way to represent the information encoded by a
higher-categorical section $\cC \to L\cE_f$. Returning to the case of $\DVect_k^\otimes \to \Fin$
with inert maps denoted as $In$, the homotopical category $\DSect_{In}(\Fin,\DVect_k^\otimes)$
provides a strictification of the infinity-category of $\rE_\infty$-algebras in $\DVect_k$ and brings with
it an alternative to \cite{LU}'s proof of the existence of (homotopy) colimits. The same
is true for more general operation-indexing categories $\rO \to \Fin$, whether or not they have
corresponding operadic analogues.

The procedure of associating $\bfE \to \bC$ to the opfibration $\cE \to \cC$ can be also
carried out in the situation when $\cE$ is an infinity-category (and the functor to $\cC$ is a
cocartesian fibration), and that is explained in subsection \ref{highercatsegsections}.
One then has that the infinity-category $\Sect(\cC,\cE)$ is identified with a full subcategory
of $\Sect(\bC,\bfE)$ in a way similar to that explained in Proposition \ref{introcomparisonbasic}.
This allows to establish the presentability part of Theorem \ref{introthmcomparisonmodel}
more generally, as stated in Theorem \ref{thmpresentabilitygeneral} in the main text.

\begin{introthm}
\label{introthmpres}
Let $\cE \to \cB$ be a cocartesian fibration over a small infinity category $\cB$, such
that each fibre $\cE(b)$ is presentable and each transition functor is accessible.
Let $\cS$ be a subset of maps of $\cB$ such that each transition functor $f_!$ along
$f \in \cS$ preserves limits (and thus satisfies the right adjoint functor theorem).

Then the infinity-category of $\cS$-cocartesian sections $\Sect_\cS(\cB,\cE)$ is presentable, with limits
and (sufficiently large) filtered colimits calculated fibrewise.
\end{introthm}
This result does not appear in \cite[5.5]{LUHTT}, and is relatively immediate
if one adapts the Segal section perspective. It gives, again, an alternative to \cite{LU}'s proof of presentability of
the categories of algebras over infinity-operads in presentably symmetric monoidal
categories.

\subsection*{Resolutions} The final section can be viewed as an exercise in the calculus of simplicial replacements,
and involves proving a certain descent statement. Given a functor $F: \cD \to \cC$,
one can readily induce the pullback $F^*: \DSect(\cC,\cE) \to  \DSect(\cD,F^*\cE)$. If
$F$ presents $\cC$ as a higher-categorical localisation of $\cD$ with respect to $F$-isomorphisms,
denoted $F-iso$, then Theorem \ref{introthmcomparisonmodel} implies that $F^*$
is homotopically fully faithful, and its image consists of $\DSect_{F-iso}(\cD,F^*\cE)$.
When $F$ belongs to a special class of localisations (called resolutions in this paper, a condition
due to \cite[Key Lemma]{HINICH}), the same statement can be proven without passing by
the comparison with the higher-categorical sections, which reveals interesting simplicial
combinatorics on the way (and that we shall use in our subsequent work).

\subsection*{Organisation of the paper} The first section introduces the combinatorics
relevant to the procedure of replacing $\cE \to \cC$ with its ``simplicial extension''
$\bfE \to \bC$, and we proceed to define and study the first properties of derived/Segal
sections, proving Propositions \ref{introcomparisonbasic} and \ref{intropropositionbasic} in the second section.
The third section compares the Segal sections with the higher-categorical sections,
which leads to Theorem \ref{introthmcomparisonmodel}, and then gives a version of
the purely higher-categorical Segal section construction, resulting in Theorem \ref{introthmpres}.
The resolutions example is treated in the last, fourth section.

\subsection*{Acknowledgements} The author is grateful to Yonatan Harpaz, Dmitry Kaledin and Carlos Simpson
for the discussions around the content of this paper. I thank both Fondation Jacques Hadamard (FMJH) and Centre de
Math\'ematiques Laurent Schwartz for financial, and my friends for personal support. This paper uses Paul Taylor's
diagrams package.

\addtocontents{toc}{\protect\setcounter{tocdepth}{2}}

\section{Simplicial Replacements}
\label{sectionsegalsections}

\subsection{Preliminaries}
For the dictionary of Grothendieck fibrations, and a definition of a factorisation category,
in particular of a Reedy category, we invite the reader to
consult the first section of \cite{BALRMSF}. Due to the existing ambiguity we shall
often interchange \emph{op} and \emph{co}
when speaking of opfibrations: they shall also be called cocartesian fibrations,
and the maps will be either cocartesian or opcartesian. The terms fibration, cartesian fibration
will mean the classical, contravariant Grothendieck fibrations.

One way to produce new factorisation categories out of existing ones
consists of considering presheaves, interpreted as discrete opfibrations.

\begin{opr}
	\label{indexedcategorydefinition}
	Let $\cB$ be any small category.
	A \emph{$\cB$-indexed category} is a (small) discrete opfibration $\cX \to \cB^\op$.
	A morphism of $\cB$-indexed categories $\cX \to \cY$ is given by
	a cocartesian morphism of discrete opfibrations over $\cB^\op$.
\end{opr}

We denote by $\Cat(\cB)$ the category of small $\cB$-indexed categories. When $\cB = [0]$
the one-object category, $\Cat([0])$ is simply denoted as $\Cat$.

\begin{rem}
	Conventionally (as for instance in topos theory \cite{JOHNSTONE}) an indexed category is yet another name for a contravariant pseudofunctor from $\cB$ to categories. We adopt a more rigid notion, which is equivalent to a presheaf of sets over $\cB$.
\end{rem}

Let now $\cB$ be a factorisation category, with the factorisation structure given by $(\sG,\sD)$.

\begin{lemma}
	For any $\cB$-indexed category $\pi:\cX \to \cB^\op$ There exists a unique factorisation system $(\sL_\cX,\sR_\cX)$ on $\cX$ such that $\pi$ becomes a factorisation functor $(\cX, \sL_\cX,\sR_\cX) \to (\cB^\op,\sD^\op,\sG^\op)$. Moreover, each morphism of $\cB$-indexed categories becomes a factorisation functor, as well.
\end{lemma}
\proof Set $\sL_\cX := \pi^{-1} (\sD^\op)$ and $\sR_\cX := \pi^{-1} (\sG^\op)$. \endproof

\begin{opr}
	\label{factorisationcanonicallyinduceddefinition}
	We shall call the pair $(\sL_\cX,\sR_\cX)$ the factorisation system
	canonically induced from $(\cB,\sG,\sD)$.
\end{opr}

\begin{ntn}
	If the factorisation category structure on $\cB$ has a name (as, for example, the Reedy factorisation system on $\Delta$), then we shall also adopt the same name for the factorisation system on the $\cB$-indexed categories $\cX \to \cB^\op$.
\end{ntn}

\begin{opr}
	Let $F: \cB' \to \cB$ be a functor. A \emph{$F$-reindexing} of a
	$\cB$-indexed category $\pi:\cX_\cB \to \cB^\op$ is the pull-back of
	$\pi$ along $F^\op$. In other words, it is the left vertical arrow in
	the pullback square
	\begin{diagram}[small]
	\cX_{\cB'} \SEpbk & \rTo^{F_\cX} 		& \cX_{\cB}	\\
	\dTo	<{\pi'}			&								&		\dTo>\pi		\\
	\cB'^\op & 		\rTo_{F^\op}		& \cB^\op	\\
	\end{diagram}
\end{opr}

\begin{lemma}
	Let $F: (\cB',\sG',\sD') \to (\cB,\sG,\sD)$ be a factorisation functor and $\cX \to \cB^\op$ a $\cB$-indexed category.  Then the functor $F_\cX:\cX_{\cB'} \to \cX_{\cB}$ induced by the reindexing operation is a factorisation functor $(\cX_{\cB'}, \sL_{\cX_{\cB'}},\sR_{\cX_{\cB'}}) \to (\cX_{\cB}, \sL_{\cX_{\cB}},\sR_{\cX_{\cB}})$ between the canonically induced factorisation systems.
\end{lemma}
\proof Immediate. \endproof

As we see, $\cB$-indexed categories naturally inherit the factorisation structure from $\cB$, and the interaction with factorisation functors is equally natural. %

\begin{prop}[Inheritance for indexed categories]
	\label{indexedcatinheritance}
	Let $F: (\cB',\sG',\sD') \to (\cB,\sG,\sD)$ be a factorisation functor. Then, for any $\cB$-indexed category $\cX_\cB$, we have the following:
	\begin{enumerate}
	\item Given any square
	\begin{diagram}[small]
	\cY_{\cB'} & \rTo^{G } 		& \cX_{\cB}	\\
	\dTo	<{\pi'}			&								&		\dTo>\pi		\\
	\cB'^\op & 		\rTo^{F^\op}		& \cB^\op	\\
	\end{diagram}
	with $\cY_{\cB'} \to \cB'^\op$ a $\cB'$-indexed category, the functor $G$ is naturally
	a factorisation functor $G:(\cY_{\cB'}, \sL_{\cY_{\cB'}},\sR_{\cY_{\cB'}}) \to (\cX_{\cB}, \sL_{\cX_{\cB}},\sR_{\cX_{\cB}})$ between the canonically induced factorisation systems.
		\item If $\sD^\op$ is a locally Noetherian category (see \cite{BALRMSF}), then so is the induced category $\sL_{\cX_{\cB}}$. There is also a dual result for the right class.
		\item If $F$ is such that the induced functor $\sD'^\op \to \sD^\op$ is a closed immersion of Noether categories (see \cite{BALRMSF}), then the induced functor $\sL_{\cX_{\cB'}} \to \sL_{\cX_{\cB}}$ has the same property, as well.
		\item If $F^\op: \cB'^\op \to \cB^\op$ is right-closed (see \cite{BALRMSF}), then so is $F_\cX: \cX_{\cB'} \to \cX_{\cB}$. Dually for left-closed.
		\item If $ (\cB,\sG,\sD)$ is a Reedy category, then so is $(\cX_{\cB}, \sL_{\cX_{\cB}},\sR_{\cX_{\cB}})$.
	\end{enumerate}
\end{prop}
\proof Most proofs are elementary, using the fact that an indexed category is discretely opfibred over the base. This for example leads to isomorphisms of comma categories $\cX_\cB / x \cong (\pi(x) \backslash \cB)^\op$, that allows to verify both the Noether property and the closed immersion condition. \endproof

\subsection{The replacements}
\begin{opr}
	\label{simplicialreplacement}
	\label{simplicialreplacementdefinition} Given a small category $\cC$, its \emph{simplicial replacement} is the unique $\Delta$-indexed category $\bC \to \Delta^\op$ such that the fibre $\bC([n])$ is the set $\Ob \Fun([n],\cC)$ of functors from $[n]$ to $\cC$, with morphisms over $[n] \leftarrow [m]$ given by precomposition $\Fun([n],\cC) \to \Fun([m],\cC)$.
\end{opr}

\begin{lemma} For $\cC \in \Cat$, the simplicial replacement $\bC \to \Delta^\op$ can be obtained as the opfibrational Grothendieck construction $\int N \cC$ of the nerve $N \cC: \Delta^\op \to \Set \subset \Cat$. The assignment $\cC \mapsto \bC$ defines a functor from $\Cat$ to the category $\Cat(\Delta)$ of $\Delta$-indexed categories. \end{lemma}
\proof Clear. \endproof

\begin{rem}
The assignment $\cC \mapsto \bC$ does not commute with the operation of taking
opposite categories. For any $\cC$, one can construct an isomorphism of categories
$\int N \cC \cong \int N (\cC^\op)$ that commutes with the $\Delta$-indexing up to
an involution reversing the order. Taking the simplicial replacements does not
send equivalences of categories to equivalences over $\Delta^\op$.
\end{rem}

\begin{ntn}
	\label{joinobject}
	An object of $\bC$ is given by a sequence $c_0 \to ... \to c_n$ of composable morphisms in $\cC$.
	It will often be denoted as $\bc_{[n]}$ or simply as $\bc$ when the $\Delta$-index is not important.
	We shall also write $\bc_{[n]}: [n] \to \cC$ for the associated functor.

	For a functor $F: \cD \to \cC$ the induced functor will often be denoted $\bF:\bD \to \bC$. One has $\bF(d_0 \to ... \to d_n)= Fd_0 \to ... \to Fd_n$.

	%
\end{ntn}

The following result already appears in \cite{TV,BALRMSF} and will be revisited in
the section devoted to resolutions.

\begin{lemma}
	\label{headandtail}
	The assignments $\bc_{[n]} \mapsto c_0$ and $\bc_{[n]} \mapsto  c_n$ determine functors
 $h_\cC:\bC \to \cC$ and $t_\cC: \bC \to \cC^\op$ that are higher-categorical localisations
 of $\bC$ along the maps that are sent by these functors to identities.
\end{lemma}
\proof The opposite of the functor $t_\cC$ is seen to be an opfibration with contractible fibres,
one then applies \cite[Lemma 3.53]{BALRMSF}. The functor $h_\cC$ is seen to be isomorphic, in $\Cat / \cC$,
to the functor $t_{\cC^\op}$, via the isomorphism relating the simplicial replacements of $\cC$ and $\cC^\op$.
\endproof

Recall (for example, from Appendix of \cite{BALRMSF}) that $\Delta$ possesses the Reedy factorisation
system given by surjective and injective maps. By Proposition \ref{indexedcatinheritance}, we get the induced
Reedy factorisation system $(\bC_{-},\bC_{+})$ on $\bC$. We shall use the term
degeneracies to call the maps belonging to $\bC_+$: those are the maps that lie over
surjections in $\Delta$.

There is another factorisation system on $\Delta$, that is given by initial element preserving
maps and left interval inclusions. Let us spell out the induced structure on $\bC$ in detail.
\begin{opr}
A map $\zeta:\bc_{[n]} \to \bc'_{[m]}$ is \emph{Segal} iff its projection in $\Delta$, $\pi(\zeta):[m] \to [n]$, is an interval inclusion of $[m]$ as first $m+1$ elements of $[n]$, i.e. $\pi(\zeta)(i) = i$ for $0 \leq i \leq m$. In particular, $m$ should be less or equal to $n$.

A map $\zeta:\bc_{[n]} \to \bc'_{[m]}$ is \emph{anchor}, or \emph{endpoint-preserving}, iff the underlying map in $\Delta^\op$ preserves the endpoints: $\pi(\zeta)(m) = n$.

\end{opr}
We denote by $\sS_\bC$ and $\sA_\bC$ the sets of all Segal and anchor maps respectively.
Note that $h_\cC$ sends $\bC_{+}$ and $\sS_\bC$ to the identity maps of $\cC$,
and $t_\cC$ sends $\bC_{+}$ and $\sA_\bC$ to the identity maps in $\cC$.
The importance of Segal maps is summarised by the following:
\begin{prop}
	\label{locprop}
	For a small category $\cC$, the functor $h_\cC: \bC \to \cC$ is a higher-categorical
	localisation of $\bC$ along $\sS_\bC$.
\end{prop}

\proof Lemma \ref{headandtail} establishes that $h_\cC$ is a localisation along the
maps of $\bC$ that project to the identities of $\cC$. It will suffice to show that
given an infinity-category $\cX$ and any functor $F: \bC \to \cX$ inverting $\sS_\bC$,
one has that $F$ also inverts all $h_\cC$-identities.

Let $\alpha : \bc\to \bc'$ be a map in $\bC$ such that $h_\cC (\alpha) = id_{c}$.
This implies that $\alpha$ takes the form
$$
\alpha: (c \to c'  \to ... \to c  \to c_0 \to ...) \longrightarrow (c \to c'_0 \to ...)
$$
with $c \to c' \to ... \to c$ composing into $c \stackrel = \to  c$. Since Segal maps are inverted by $F$, it is clear that the first element preserving maps are $F$-invertible.
It will thus suffice to show that the following map
in $\bC$ is $F$-invertible:
$$
\alpha': (c \stackrel = \to  c  \to c_0 \to ...) \longrightarrow (c \to c'_0 \to ...),
$$
as there is a pair of endpoint-preserving maps connecting $\alpha$ and $\alpha'$ by the means of a commutative diagram.
In turn it will suffice to show that the following map (connected to $\alpha'$ by Segal maps) in $\bC$ is $F$-invertible:
$$
\beta: (c \stackrel = \to c) \longrightarrow (c).
$$
The map $\beta$ is a one-sided inverse of the degeneracy $(c) \to (c \stackrel =  \to  c)$.
The latter is a one-sided inverse of the Segal map $(c \stackrel =  \to c) \to (c)$. All these
maps are hence $F$-invertible.
\endproof
\begin{rem}
 Neither the set of Segal maps nor the set of $h_\cC$-identities is saturated
 in the sense one applies when one speaks of localisation \cite{DHKS}: the saturated
 set of maps that gives the localisation $h_\cC: \bC \to \cC$ is given by those
 $\alpha:\bc_{[n]} \to \bc'_{[m]}$ that project to isomorphisms of $\cC$.
\end{rem}

Proposition \ref{locprop} justifies the idea that, given a homotopical category $(\cM,\cW)$, a functor $F: \bC \to \cM$ sending $\sS_\bC$ to $\cW$ is a suitable weakening of the concept of a functor from $\cC$ to $\cM$.
The action of $F$ on spans in $\bC$ like
$$
c  \longleftarrow (c \stackrel f \to c' ) \longrightarrow c',
$$
where the left arrow is Segal, gives a span $F(c) \stackrel \cW \leftarrow F(c \to c') \rightarrow F(c')$, where the left map is a weak equivalence. On the level of $\Ho \cM$, this span gives a map $F(c) \to F(c')$, which one can denote $F(f)$. Applying $F$ to higher-length objects then ensures higher coherences for the `weak functor' $F$.
This is summarised by saying that given a relative \cite{BARKAN} functor $F:(\bC,\sS_\bC) \to (\cM,\cW)$ (meaning that $F$ takes $\sS_\bC$ to $\cW$),
after passing to infinity-localisations one gets an induced infinity-functor $LF: \cC \cong L_{\sS_\bC} \bC \to L_\cW \cM$.

\subsection{Families over the simplicial replacement}

We now turn to the question of extending an opfibration $\cE \to \cC$ to $\bC$.
This can be done in a few ways. Since we have a functor $h_\cC: \bC \to \cC$,
there is a pull-back opfibration $h_\cC^* \cE \to \bC$. However, unless $\cE \to \cC$
is a Quillen presheaf, it is not apparent why the category of sections of $h_\cC^* \cE \to \bC$
carries a model structure. In particular, we cannot use this approach to study higher algebra
in the language of model categories.

Any opfibration $p:\cE \to \cC$ can be described, up to an
equivalence, as a Grothendieck construction of a covariant functor from $\cC$ to categories.
The latter is the same thing as a contravariant functor
from $\cC^\op$ to categories. The way to capture this duality without
passing to category-valued functors, is the following.

\begin{opr}
	\label{transposefib}
	Fix an opfibration $p: \cE \to \cC$. Define a category denoted as $\cE^\top$ as follows:
	\begin{enumerate}
		\item $Ob(\cE^\top) = Ob(\cE)$
		\item A morphism from $x \to z$ in $\cE^\top$ is an isomorphism class of cospans in $\cE$
		$$
		x \longrightarrow y \longleftarrow z
		$$
		such that the left arrow is fibrewise, $p(x \to y) = id_{p(x)}$, and the right arrow is opcartesian.
	\end{enumerate}
	There is an evident functor $p^\top: \cE^\top \to \cC^\op$ which sends a map
  $x \longrightarrow y \longleftarrow z$ to $p(y \longleftarrow z)$.
  A morphism of $\cE^\top$ is $p^\top$-cartesian iff it can be represented
  by a span of the form $y \stackrel{=}{\longrightarrow} y \longleftarrow z$.
  The functor $p^\top$ is a fibration, which we call
  the \emph{transpose fibration} of $p$.
\end{opr}

If $\cE \to \cC$ equals $\int E \to \cC$ for a functor $E: \cC \to \Cat$, then $\cE^\top \to \cC^\op$ is equivalent to the (fibrational) Grothendieck construction applied to $E: (\cC^\op)^\op \to \Cat$ viewed as a contravariant functor on $\cC^\op$. In particular, if $\cE \cong \cE_0 \times \cC$ and the opfibration
structure is given by the $\cC$-projection, then $\cE^\top$ is simply $\cE_0 \times \cC^\op$.
The construction of a transpose fibration that we outlined can be also done higher-categorically,
see \cite{BARGLANAR}. 


After taking the transpose fibration $\cE^\top \to \cC$ of the opfibration $\cE \to \cC$,
we can also consider the pull-back fibration $t_\cC^* \cE^\top \to \bC$. The following
proposition is a manifestation of the duality between $\cE$ and $\cE^\top$:

\begin{prop}
	\label{EtoEtop}
	Given an opfibration $p: \cE \to \cC$, there is a morphism $T:h^*_\cC \cE \to t^*_\cC \cE^\top$ commuting with functors to $\bC$ which sends opcartesian maps of $h^*_\cC \cE$ to cartesian maps of $t^*_\cC \cE$. 
\end{prop}
\proof Consider the category $\cX$ defined as follows.
\begin{itemize}
	\item An object of $\cX$ is a pair $(\bc_{[n]}, \alpha)$ where $\bc_{[n]} = c_0 \to ... \to c_n$ is an object of $\bC$ and $\alpha: x \to y$ is an opcartesian map in $\cE$ which covers the composition $c_0 \to c_n$ in $\cC$ (i.e. $p(\alpha) = c_0 \to c_n$),
	\item A morphism $(\bc_{[n]}, \alpha:x \to y) \to (\bc'_{[m]}, \beta: x' \to y')$ consists of a map $\bc \to \bc'$ in $\bC$ and a map $\gamma: x \to x'$ which covers the induced map $c_0 \to c'_0$.
\end{itemize}
One can check that the natural functor $\cX \to \bC$ is an opfibration, and that the assignment $(\bc, \alpha: x \to y) \mapsto (\bc, x)$ defines an equivalence over $\bC$ of opfibrations $\cX \stackrel \sim \to h_\cC^* \cE$.

On the other hand, consider the assignment $(\bc, \alpha: x \to y) \mapsto (\bc, y)$. We claim that it defines a functor $\bar T: \cX \to t^*_\cC \cE$ commuting with projections to $\bC$. Let $(f,t):(\bc, \alpha: x \to y) \to (\bc', \beta: x' \to y')$ be a map. In particular, we have the following diagram in $\cE$:
\begin{equation}
\label{diagramfort}
\begin{diagram}[small,nohug]
x 						& \rTo^t 	& x'	 				\\
\dTo<\alpha		&				&	\dTo>\beta \\
y						&				& 	y'.		\\
\end{diagram}
\end{equation}
Suppose first that the map $t$ is fibrewise. Then by opcartesian property there exists a map $t': y \to y'$ rendering the diagram (\ref{diagramfort}) commutative. Remembering the description of arrows in Definition \ref{transposefib}, we define $\bar T(f,t) = (f,y \stackrel{t'}{\to} y' \stackrel{id}{\leftarrow} y')$; in other words, we view $t'$ as a fibrewise map of $\cE^\top$.

Next, if $t$ is opcartesian, find an opcartesian map $k: y' \to z$ in $\cE$ covering $c'_m \to c_n$ (which is induced from $f:\bc \to \bc'$). The composition $k \beta t$ and $\alpha$ both project along $\cE \to \cC$ to the map $c_0 \to c_n = c_0 \to c'_0 \to c'_m \to c_n$, hence there is a (fibrewise) isomorphism $z \cong y$. This implies that the diagram (\ref{diagramfort}) can be completed as
\begin{diagram}[small,nohug]
x 						& \rTo^t 	& x'	 				\\
\dTo<\alpha		&				&	\dTo>\beta \\
y						&	\lTo^{t'}			& 	y'.	\\
\end{diagram}
with all arrows opcartesian in $\cE$. We put, again, $\bar T(f,t) = (f,y \stackrel{id}{\to} y \stackrel{t'}{\leftarrow} y')$, thus viewing $t'$ as a cartesian map of $\cE^\top$. Any other case of $(f,t)$ can be treated by reducing to these two cases.

Inverting the equivalence $\cX \stackrel \sim \to h_\cC^* \cE$ and composing with $\bar T$, we obtain the desired functor $T: h_\cC^* \cE \to t^*_\cC \cE$.
\endproof

As we can see from the construction, the functor $T: h_\cC^* \cE \to t^*_\cC \cE$ acts by pushing an object $x$ over $c_0$ along the composition of the string $$c_0 \stackrel{f_1}{\to} ... \stackrel{f_n}{\to} c_n,$$ which gives $(f_n ... f_1 )_! x \in \cE(c_n)$. There exists a refinement of this procedure that remembers much more information.

\begin{opr}
Given $\bc_{[n]} \in \bC$, consider it as a functor $\bc_{[n]}^\op : [n]^\op \to \cC^\op$,
and define $\bfE(\bc_{[n]})$ to be $\Sect([n]^\op, \bc_{[n]}^{\op,*} \cE^\top)$, the category of
sections of the transpose fibration pulled back along $\bc_{[n]}^\op$.
\end{opr}

\begin{opr}
	\label{simplicialextensionofanopfibration}

	For an opfibration $\cE \to \cC$, its \emph{simplicial extension} is an opfibration $\bfE \to \bC$ which is the (covariant) Grothendieck construction of the functor $\bc \mapsto \bfE(\bc)$. It is characterised by the fibres $\bfE(\bc_{[n]}) \cong \Sect([n]^\op,  \bc_{[n]}^{\op,*}  \cE^\top)$ and by the transition adjunctions $\alpha_!: \bfE(\bc_{[n]}) \rightarrow \bfE(\bc'_{[m]})$ given by pullbacks along maps $\alpha: \bc_{[n]} \to \bc'_{[m]}$.
\end{opr}

\begin{lemma}
	\label{simplextensionprelim} %

	The assignment $\bc_{[n]} \mapsto \bfE(\bc_{[n]})= \Sect([n]^\op, \bc_{[n]}^{\op,*} \cE^\top)$ has the following properties:
	for each Segal map $\alpha: \bc_{[n]} \to \bc'_{[m]}$ there is an induced adjunction
	$$
	\alpha_!: \bfE(\bc_{[n]}) \rightleftarrows \bfE(\bc'_{[m]}) : \alpha^*.
	$$
 When the fibres of $\cE \to \cC$ are complete, the right adjoint $\alpha^*$ exists for
 all maps $\alpha$ of $\bC$.
\end{lemma}
\proof
Given $\alpha: \bc_{[n]} \to \bc'_{[m]}$, the functor $	\alpha_!: \bfE(\bc_{[n]}) \rightarrow \bfE(\bc'_{[m]})$ is induced as a pullback along the underlying map $[m] \to [n]$. It will thus suffice to prove that each $\alpha_!$ has an adjoint. To do this, consider four different cases:
\begin{enumerate}
	\item The map $\alpha$ is a degeneracy, covering a surjection $a:[m] \twoheadrightarrow [n]$. Given a section $A \in \bfE(\bc'_{[m]})$, which can be represented by a diagram $A_m \to ... \to A_0$ (we use this order to underline the opposite character), the section $\alpha^* A$ is then represented by the diagram
	$$
	A_{a^{-1}(m)} \to ... \to A_{a^{-1}(0)}
	$$
	where $a^{-1}(i)$ denotes the first element of the inverse image of $i \in [n]$ under $a$.


	\item The map $\alpha$ covers an injection $a: [m] \hookrightarrow [n]$ which preserves both endpoints. We use $a$ to identify $[m]$ with its image in $[n]$. Representing again
  $A \in \bfE(\bc'_{[m]})$ as $A_m \to ... \to A_0$, we set
	where $A(j) \to B_i \to A(k)$ is the factorisation of the arrow $A(j) \to A(k)$, with $j$ being the maximal closest to $i$ (respectively $k$ being the minimal closest to $i$) value belonging to $[m]$. Intuitively, the procedure consists of ``filling the holes'', the elements which are not in the image of $a$, in the natural and minimal way.

  \item For the remaining two injections, the one omitting the first or the last
  element respectively, the situation is as follows. The Segal maps
	$$
  (c_0 \to ... \to c_m \stackrel f \to c_{m+1}) \, \, \, \longrightarrow (c_1 \to ... \to c_m)
  $$
	omitting
  the last element are treated similarly to the previous case: from
  $$
  A_0 \longleftarrow ... \longleftarrow A_{m}
  $$
  we get a section
  $$
  A_0 \longleftarrow ... \longleftarrow A_{m} \longleftarrow f_! A_{m}.
  $$
  Suppose we have a map in $\bC$ that looks like
  $$
  (c_0 \to ... \to c_n) \, \, \, \longrightarrow (c_1 \to ... \to c_n);
  $$
  given $A_1 \leftarrow ... \leftarrow A_n$, we define $B_0 \leftarrow ...
  \leftarrow B_n$ by setting $B_0=*$ to be the terminal object of $\cE(c_0)$.
  Next, $B_1 := A_1 \times f_{1,!} *$, to correct for the lack of map to
  $f_{1,!} *$. The procedure then continues, with
  $$
  B_2 := A_2 \prod_{f_{2,!}A_1} f_{2,!} B_1 \text{ and }
  B_k := A_k \prod_{f_{k,!}A_{k-1}} f_{k,!} B_{k-1},
  $$
  which is indeed the application of inductive principle to construct sections,
  providing us with the right adjoint just like in
	Proposition \cite[Proposition 1.31]{BALRMSF}. Intuitively, the fact that the
	terminal object is not preserved by transition functors forces us to ``correct''
	the rest of the section. Due to the fact that this rarely happens in practice,
	the procedure is in fact trivial, and consists of putting the final object in
	$c_0$-th position. \endproof
\end{enumerate}


The functor $\bfE \to \bC$ will be our chosen way to extend $\cE \to \cC$ to
$\bC$. It is an opfibration and always a fibration over the Segal maps (and a bifibration
if $\cE \to \cC$ is fibrewise-complete).

\begin{rem}
	Note that $\bfE$ is not a simplicial replacement of $\cE$ or $\cE^\top$. In particular, the fibre of $\bfE \to \bC$ over an object $\bc_{[n]}$ is equivalent to $\Sect([n]^\op, \bc_{[n]}^* \cE^\top)$, with $\bc_{[n]}$ regarded as a functor $[n]^\op \to \cC^\op$.
	One could view such sections over an $n$-simplex as a relativisation and a thickening of the nerve construction. For example, when $\cE \cong \cE_0 \times \cC$ is a constant opfibration, then $\bfE(\bc_{[n]})$ is simply the category $\Fun([n]^\op,\cE_0)$.

\end{rem}

Given any $\bc_{[n]} \in \bC$, the map $\bc_{[n]} \to c_0$ is Segal. In the induced adjunction,
\begin{equation}
\label{remadjseg}
\bfE(\bc_{[n]}) \rightleftarrows \bfE(c_0) = \cE(c_0),
\end{equation}
observe that the right adjoint takes values in
the full subcategory $\bfE^{cart}(\bc_{[n]})=\Cart([n]^\op,\bc_{[n]}^* \cE^\top)$ spanned by cartesian sections of $\Sect([n]^\op,\bc_{[n]}^* \cE^\top)$.
Note that for any map $\alpha:\bc \to \bc'$, the left adjoint preserves cartesian sections, giving
$\alpha_! :\bfE^{cart}(\bc) \to \bfE^{cart}(\bc')$. The induced opfibration $\bfE^{cart} \to \bC$ is in
fact equivalent to $h_\cC ^* \cE \to \bC$.


\begin{lemma}
	\label{simpreplacementrelation}
	Let $\cE \to \cC$ be an opfibration. Then there is a diagram $h_\cC^* \cE \stackrel S \to \bfE \to t_\cC^* \cE^\top$  of functors over $\bC$ factoring the functor $T$ of Proposition 	\ref{EtoEtop}.
	The functor $S$ preserves opcartesian arrows, and for each $\bc$, the induced functor
	$S_\bc: \cE(c_0) \cong h_\cC^* \cE(\bc) \to \bfE(\bc)$ is fully faithful. The image of $S$
	consists of the \emph{subopfibration} $\bfE^{cart} \subset \bfE$ given by the subcategories
	of cartesian sections in $\bfE(\bc_{[n]})$.
\end{lemma}
\proof  The functors of the second assertion are defined as follows. Given $\bc_{[n]}: [n]^\op \to \cC^\op$, the second functor $\bfE \to t_\cC^* \cE^\top$ is induced by evaluating a section $X \in \Sect([n]^\op,\bc_{[n]}^* \cE^\top)$ at $n$, which corresponds to $c_n$. The first functor, $S: h_\cC^* \cE \to \bfE$,
is induced by the inclusion $\bfE^{cart} \subset \bfE$ and the equivalence $\bfE^{cart} \cong h_\cC^* \cE$.
It can be constructed explicitly via the procedure that sends $X \in \cE(c_0)$ to the following diagram in $\cE^\top$:
$$
X \longleftarrow (f_1)_! X \longleftarrow ... \longleftarrow (f_n...f_1)_! X;
$$
here $\bc_{[n]} = c_0 \stackrel{f_1}{\longrightarrow} ... \stackrel{f_n}{\longrightarrow} c_n$ and $(f_i)_!$ and the like denote the transition functors of $\cE \to \cC$. The verification of the rest is elementary. \endproof



\begin{rem}
Given an opfibration $p:\cE \to \cC$ and a string of composable arrows $\bc_{[n]}:[n] \to \cC$,
one can consider the category of sections $\Sect([n],\bc_{[n]}^* \cE)$ directly induced from
the opfibration $p$. Without going into much detail, we note that the bifibration over $\int N(\cC^\op)$
resulting from the assignment $\bc_{n} \mapsto \Sect([n],\bc_{[n]}^* \cE)$ is related to studying (weak) sections
of the transpose fibration $\cE^\top \to \cC^\op$.
\end{rem}

\section{Derived, or Segal, sections}

\subsection{Presections}

\begin{opr}
	Given an opfibration $\cE \to \cC$, its category of \emph{presections} is the category
	$$
	\PSect(\cC,\cE):= \Sect(\bC, \bfE).
	$$
	of sections of the simplicial extension $\bfE \to \bC$.
\end{opr}

To relate $\Sect(\cC,\cE)$ to $\PSect(\cC,\cE)$, recall the functor $S:h_\cC^* \cE \to \bfE$ of Lemma \ref{simpreplacementrelation}. The functor $h_\cC$ also induces the pull-back functor $h_\cC^*: \Sect(\cC,\cE) \to \Sect(\bC,h_\cC^* \cE)$.

\begin{prop}
	\label{sectinpsect}
	The assignment $X \mapsto S \circ (h_\cC^* X)$ defines a functor $i:\Sect(\cC,\cE) \to \PSect(\cC,\cE)$. Its essential image consists of the presections sending the Segal maps $\sS_\bC$ to \emph{cartesian} morphisms in $\bfE$.
\end{prop}
\proof
The functor $h_\cC: \bC \to \cC$ is a localisation along Segal maps. This readily implies that
the natural pull-back functor $\Sect(\cC,\cE) \to \Sect(\bC,h_\cC^*\cE)$ is fully faithful,
and its image consists of sections $X: \bC \to h_\cC^* \cE$ that send the Segal maps
of $\bC$ to the opcartesian maps of $h_\cC^* \cE$. Lemma \ref{simpreplacementrelation} then
implies that the category $\Sect(\cC,\cE)$ is identified with a full subcategory of
$\Sect(\bC,\bfE)$ consisting of such sections $X$ that
\begin{enumerate}
\item the section $X$ factors through $\bfE^{cart} \to \bC$: for each $\bc \in \bC$, the value $X(\bc)$ is a cartesian section in $\bfE(\bc)$,
\item for each Segal map $\bc \to \bc'$, the induced map $X(\bc) \to X(\bc')$ is opcartesian.
\end{enumerate}
In the presence of the first condition, the second condition is equivalent to requiring
that for each Segal map $\bc \to \bc'$, the induced map $X(\bc) \to X(\bc')$ is \emph{cartesian}.
This happens because the adjunction $\bfE(\bc) \rightleftarrows \bfE(\bc')$ restricts
to an adjoint equivalence of categories $\bfE(\bc)^{cart} \cong \bfE(\bc')^{cart}$ whenever
the map $\bc \to \bc'$ is Segal. Thus for Segal $\alpha:\bc \to \bc'$, the map
$ \alpha_! X(\bc) \to X(\bc')$ is an isomorphism iff $  X(\bc) \to \alpha^* X(\bc')$ is such.

However, any section $X:\bC \to \bfE$ sending
the Segal maps of $\bC$ to cartesian maps of $\bfE$ automatically takes values in $\bfE^{cart}$,
as verified by computing the value of $X$ on $\bc \to c_0$.
\endproof
\begin{rem}
	Consider an object $\bc_{[n]} = c_0 \stackrel{f_1}{\longrightarrow} c_1 \stackrel{f_2}{\longrightarrow} ... \stackrel{f_n}{\longrightarrow} c_n$ of $\bC$. Then $S \in Sect(\cC, \cE)$ is sent by the functor above to $i(S)$ such that $i(S)(\bc_{[n]})$ is represented by the diagram
	$$
	S(c_0) \longleftarrow (f_1)_! S(c_0) \longleftarrow ... \longleftarrow (f_n ... f_1)_! S(c_0)
	$$
	where $(f_k... f_1)_!: \cE(c_0) \to \cE(c_k)$ is a transition functor along the composition of $f_i$.
\end{rem}

We now put a homotopical structure on $\cE \to \cC$. The notion that leads to the
model structure on presections is the following one:
\begin{opr}
	\label{modelfibrationdefinition}
	A \emph{model opfibration} $\cE \to \cC$ is an opfibration such that each fibre $\cE(c)$ is a model category and the transition functors preserve fibrations and trivial fibrations of the model structure.
	Equivalently, given a diagram
	\begin{diagram}[small]
	X			& \rTo	& Y \\
	\dTo		&			& \dTo \\
	Z			& \rTo	& T
	\end{diagram}
	with horizontal maps opcartesian and vertical maps in fibres, if $X \to Z$ is
	a (trivial) fibration then so is $Y \to T$.
\end{opr}

\begin{opr}
\label{defantisegalandconvex}
Say that a map $f:\bc_{n} \to \bc_{m}$ is
\begin{enumerate}
\item anti-Segal if the underlying $\Delta$-map is an interval inclusion of $[m]$
as last $m+1$ elements of $[n]$,
\item convex if the underlying $\Delta$-map preserves initial and final elements.
\end{enumerate}
\end{opr}

\begin{lemma}
Any map $\alpha: \bc \to \bc'$ can be factored as $\bc \to \bc'' \to \bc'$, where the first map is anti-Segal and the second map preserves the initial elements.
\end{lemma}
\proof Clear. \endproof

\begin{lemma}
\label{simplextisasegalbifibr}
\label{lmsegalbasechage}
Let $\cE \to \cC$ be a model opfibration. Then
\begin{enumerate}
	\item The bifibration $\bfE \to \bC$ is a Quillen presheaf \cite{H-S}:
	the fibres $\bfE(\bc_{[n]})=\Sect([n]^\op, \bc_{[n]}^{\op,*} \cE^\top)$ are model categories with cofibrations and weak equivalences given
	valuewise, and for each $f: \bc \to \bc'$,
	the adjunction $f_!: \bfE(\bc) \rightleftarrows \bfE(\bc'): f^*$ is a Quillen pair.
	\item The following base change condition holds. For a commutative square in
	$\bC$,
	\begin{diagram}[small]
	x & \rTo^f & y \\
	\dTo<\alpha & &  \dTo>\beta \\
	z & \rTo_g & t, \\
	\end{diagram}
with vertical arrows Segal, and horizontal arrows projecting to initial
element preserving maps in $\Delta$, the induced derived natural transformation
$\bL f_! \bR \alpha^* \to \bR \beta^* \bL g_!$ is an isomorphism.
	\item The following base change condition holds. For a commutative square in
	$\bC$,
	\begin{diagram}[small]
	x & \rTo^f & y \\
	\dTo<\alpha & &  \dTo>\beta \\
	z & \rTo_g & t, \\
	\end{diagram}
	with vertical arrows anti-Segal, and horizontal arrows convex, the induced
	derived natural transformation $\bL \alpha_! \bR f^* \to \bR g^* \bL \beta_!$
	is an isomorphism.
\end{enumerate}
\end{lemma}
\proof The first assertion is an application of \cite[Theorem 1]{BALRMSF}
to the fibres of $\bfE$. Observe that the fibres $\bfE(\bc_{[n]})$ are the
categories of sections of the fibration $\bc^*_{[n]} \cE^\top \to [n]^\op$,
which is also an admissible model semifibration in the terminology of \cite{BALRMSF}. The admissibility is true because
the category $[n]^\op$ and the associated matching categories have initial objects.
We thus can apply the aforementioned theorem
to establish the existence of a model structure on $\bfE(\bc)$ in which cofibrations and weak
equivalences are defined objectwise. This readily implies that the functors
$$f_!: \bfE(\bc) \longrightarrow \bfE(\bc')$$ are left Quillen.

The second and the third assertion follow from the explicit description
of functors in Lemma \ref{simplextensionprelim}. \endproof

\begin{corr}
	\label{modelstructureonpresections}
	Let $\cE \to \cC$ be a model opfibration, then the presection category $\PSect(\cC,\cE) = \Sect(\bC,\bfE)$ has the Reedy model structure of \cite[Theorem 1]{BALRMSF}.  \endproof
\end{corr}

\begin{rem}
\label{remexamples}
The fact that in Definition \ref{modelfibrationdefinition}, the transition functors preserve the fibrational part, and not
the cofibrations, reflects the dual, ``coalgebraic'' character of the notions of a presection and of a Segal
section (introduced below).

If we consider opfibrations $\cM^\otimes \to \Fin$
associated to monoidal categories as explained in the introduction, then to get a
model opfibration, one can take $\cM$ to be a Cartesian monoidal model category like
$\SSet$, but also the categories $\DVect_k$ of chain complexes of vector spaces, or simplicial
vector spaces, since over a field $k$ the tensor product is exact (note that there is no restriction
imposed on the characteristic). Another example involves
taking $\cM = \cN^\op$, where $\cN$ is a monoidal model category with all objects cofibrant. For
example, we can take $\cM = \DVect_k^\op$.
\end{rem}

\subsection{Homotopical category of Segal sections}
\label{homcatofsegsect}





\begin{opr}
\label{defsegalsection}
	Let $\cE \to \cC$ be a model opfibration. A presection $S: \bC \to \bfE$ is
	\emph{derived}, or \emph{Segal}, if for any Segal map $\alpha:\bc \to \bc'$ the
	morphism $S(\bc) \to S(\bc')$ is weakly cartesian, meaning that
the	induced morphism $S(\bc) \to \bR \alpha^* S(\bc')$ is an isomorphism in $\Ho \bfE(\bc)$.

\end{opr}

We denote by $\DSect(\cC,\cE)$ the full subcategory of $\PSect(\cC,\cE)$
consisting of Segal sections, with the weak equivalences being those induced from presections.

\begin{lemma}
	\label{equivalenttosegalissegal}
	Let $S \to S'$ be a weak equivalence in $\PSect(\cC,\cE) = \Sect(\bC,\bfE)$.
	Then, if one
	of $S,S'$ is Segal, so is the other. \endproof
\end{lemma}


\begin{lemma}
\label{derivedsectionsdontseedegenerations}
\label{lmsegalsectionsandinitialelementmaps}
Let $X$ be a derived section and $\alpha: \bc \to \bc'$ be a map in $\bC$ covering
an initial element preserving map. Then both $X(\bc) \to \bR \alpha^* X(\bc')$
and $\bL \alpha_! X(\bc) \to X(\bc')$ are isomorphisms in corresponding homotopy categories.
\end{lemma}
\proof Since $\alpha$ preserves initial elements, we have the following commutative diagram
\begin{diagram}[small,nohug]
\bc & \rTo^\alpha & \bc' \\
\dTo<p	&			&			\dTo>{q} \\
c_0 & \rTo^= & c'_0 \\
\end{diagram}
using this diagram, the base change of Lemma \ref{simplextisasegalbifibr}  and
the fact that $X$ is Segal, we
can base change from the two maps of this lemma to studying what
happens at $c_0 = c'_0$.
And there, everything trivially coincides. \endproof

\begin{rem}
Using a similar argument, it is easy to see that a presection $X$ is Segal iff
it satisfies the condition of Definition \ref{defsegalsection} for Segal maps of
the form $\bc_{[n]} \to c_0$.
\end{rem}

Consider the following span in $\bC$:
\begin{diagram}[small,nohug]
&						&	c \stackrel f \longrightarrow c'	&				&	\\
&	\ldTo^{\alpha}	&										&	\rdTo^\beta	&	\\
c 		&						&										&				& c'
\end{diagram}
A derived section $X$ would provide us with the span
\begin{diagram}[small,nohug]
&						&	\bL \beta_! X(c \stackrel f \longrightarrow c')	&				&	\\
&	\ldTo^{\sim}	&										&	\rdTo	&	\\
\bL \beta_! \bR \alpha^* X(c) 		&						&										&				& X(c')
\end{diagram}
so that on homotopy level, one has a morphism
$\bL \beta_! \bR \alpha^* X(c) \to X(c')$. However, in this case, $\bL \beta_! \bR \alpha^* X(c)$
is isomorphic to $\bR f_! X(c)$, where $\bR f_!: \Ho \cE(c) \to \Ho \cE(c')$
is the derived transition functor along $f$. That is, the data of a derived
section give us a section of the opfibration $\Ho \cE \to \cC$, with $\Ho$ applied
fibrewise:
\begin{lemma}
\label{lmhosection}
Let $X \in \DSect(\cC,\cE)$ be fibrant as a presection. Then the assignment
$c \mapsto X(c)$ naturally prolongs to a section $hX \in \Sect(\cC,\Ho(\cE))$,
where the opfibration $\Ho \cE \to \cC$ has localised categories as fibres and
right derived functors of transition functors of $\cE$, as transition functors.
\end{lemma}
\proof One can check the compositions using Lemmas \ref{lmsegalbasechage}
and \ref{lmsegalsectionsandinitialelementmaps}. \endproof

\begin{lemma}
\label{lmderivedsectionslocallyconstantalongiso}
Let $X$ be a derived section. Then for any anti-Segal map $\alpha:\bc_{[n]} \to \bc'_{[m]}$
such that the maps $c_{i-1} \to c_i$, $1 \leq i \leq n-m$, are isomorphisms, the induced map
$$
\bL \alpha_! X(\bc_{[n]}) \to X(\bc'_{[m]})
$$
is an isomorphism in $\Ho \bfE(\bc'_{[m]})$.
\end{lemma}

\proof Lemmas \ref{lmsegalbasechage} and \ref{lmsegalsectionsandinitialelementmaps} allow
to reduce everything to the case of a single isomorphism, $\bL \alpha_! X(c_0 \cong c_1) \to X(c_1)$,
and the latter is invertible in the homotopy category due to Lemma \ref{lmhosection}. \endproof

Lemma \ref{lmderivedsectionslocallyconstantalongiso} motivates to consider the following
definition.
Denote by $\cS$ a subset of maps of $\cC$.

\begin{opr}
	\label{deflocconstdersect}
	Let $\cE \to \cC$ be a model opfibration. A derived section $X \in \DSect(\cC,\cE)$ is \emph{$\cS$-locally constant} if for any anti-Segal map $\alpha:\bc_{[n]} \to \bc'_{[m]}$ such that the maps $c_{i-1} \to c_i$, $1 \leq i \leq n-m$ belong to $\cS$, the induced map
	$$
	\bL \alpha_! X(\bc_{[n]}) \to X(\bc'_{[m]})
	$$
	is an isomorphism in $\Ho \bfE(\bc'_{[m]})$.
\end{opr}

The anti-Segal maps appearing in this definition will sometimes be called $\cS$-decolouring.
We denote by $\DSect_\cS(\cC,\cE)$ the full subcategory of $\cS$-locally constant
derived sections.

\begin{lemma}
\label{lemmalocconstsecconsistency}
A derived section $X \in \DSect(\cC,\cE)$ is $\cS$-locally constant iff for any
map $\alpha:\bc \to \bc'$ projecting to an element of $\cS$, the induced
map $\bL \alpha_! X(\bc ) \to X(\bc' )$ is an isomorphism in $\Ho \bfE(\bc')$.
\end{lemma}
\proof Use of base-change similar to the proof of previous lemmas coupled with the
fact that any map of $\bC$ factors into the composition of an anti-Segal map and an
initial element preserving map.  \endproof

\begin{lemma}
\label{lmantisegmapswhenlimitsarepreserved}
Let $\cE \to \cC$ be a model opfibration and $\cS$ a subset of maps of $\cC$ such that
for each $f \in \cS$, the induced transition functor $f_!$ preserves limits.
Then for each anti-Segal map $\alpha:\bc_{[n]} \to \bc'_{[m]}$ such that the maps $c_{i-1} \to c_i$, $1 \leq i \leq n-m$ belong to $\cS$, the induced functor $\alpha_!:\bfE(\bc_{[n]}) \to \bfE(\bc'_{[m]})$
preserves fibrations, trivial fibrations and limits, and hence homotopy limits.
\end{lemma}
\proof It follows from \cite[Proposition 1.27]{BALRMSF} that the limits in $\bfE(\bc_{[n]})$
are calculated objectwise in $c_0,...c_{n-m}$ in this particular case, and after that by the means of the same inductive
procedure as in $\bfE(\bc'_{[m]})$. The observation for the (trivial) fibrations
is similar, given their explicit description \cite[Definition 2.8]{BALRMSF}. The conclusion
about the homotopy limits follows from the general results of \cite{DHKS} that reduce
the question to the computation of limits of Reedy-fibrant diagrams that are preserved by
the functors like $\alpha_!$. \endproof

\begin{prop}
	\label{segalsectionsbehaveasintended}
 Let $\cE \to \cC$ be a model opfibration, $\cS$ a subset of maps of $\cC$, and $\bfE \to \bC$ the associated
 bifibration. Then
	\begin{enumerate}
		\item if $X \in \DSect_\cS(\cC,\cE)$, then any fibrant and cofibrant replacement
		of $X$ is also a $\cS$-locally constant derived section,
		\item if $X \in \DSect_\cS(\cC,\cE)$ is fibrant as a section in
		$\PSect(\cC,\cE)=\Sect(\bC,\bfE)$ and $f:x \to y$ is a Segal map in $\bC$,
		then the induced morphism
		$X(x) \to f^* X(y)$ is a trivial fibration of fibrant objects,
		\item if $X_\bullet: I \to \DSect(\cC,\cE)$ is a diagram of derived sections,
		then its homotopy limit in $\PSect(\cC,\cE)$ is a derived section. If the transition
		functors along $\cS$ preserve limits, then the same result is true for the diagrams
		valued in $\DSect_\cS(\cC,\cE)$.
	\end{enumerate}
\end{prop}
\proof
The first assertion is a direct consequence of Lemma \ref{equivalenttosegalissegal}.

For the second assertion, we know that $X$ fibrant implies $X(x)$ fibrant for each $x \in \cX$ of degree $0$. In general, we know that $X(x) \to \Mat_x X$ is a fibration. Similarly as for simplicial objects in a model category, one can prove that for any $\alpha:x \to y$ covering an injection $[n] \hookleftarrow [m]$ in $\Delta$, the map $\Mat_x X \to \alpha^* X(y)$ is a fibration. This implies that $X(x)$ is fibrant and any Segal map $x \to y$ goes to $X(x) \to X(y)$, a fibration and a weak equivalence.

The final assertion is easily checked due to the fact that the (homotopy) limits
are calculated fibrewise \cite{H-S} in $\Sect(\bC,\bfE)$ and the commutativity between $\bR f^*$
for a Segal map $f: x \to y$ and homotopy colimits:
$$
\mathbb R \lim X_\bullet(x) \stackrel \sim \to \mathbb R \lim \mathbb R f^* X_\bullet(y) \cong
 \mathbb R f^*  \mathbb R \lim X_\bullet(y).
$$
When the transition functors along $\cS$ preserve limits, Lemma \ref{lmantisegmapswhenlimitsarepreserved}
and a dual calculation permit to check the same for the local constancy condition.
\endproof

The colimits of derived sections, if exist,
are usually not calculated in a fibrewise way. However, some abstract existence results can
be established if the original model opfibration $\cE \to \cC$ has combinatorial fibres
and accessible transition functors. For example, one can state the following:

\begin{prop}
\label{propdsectasacmc}
Let $\cE \to \cC$ be a model opfibration such that each $\cE(c)$ is left proper combinatorial,
and each transition functor $f_!: \cE(c) \to \cE(c')$ is accessible. Then
there exists a left proper combinatorial model structure on $\PSect(\cC,\cE)$ with
Reedy cofibrations and fibrant objects given by Reedy-fibrant objects of $\DSect(\cC,\cE)$.
\end{prop}
\proof According to \cite[Proposition 2.30]{BALRMSF}, each model category $\bfE(\bc)$,
being the category of sections over a Reedy category, is combinatorial. It is moreover
left proper, as is easily checked using the fact that the colimits are computed fibrewise.
Thus the Quillen presheaf $\bfE \to \bC$ has left proper combinatorial fibres. This
again implies that the category $\PSect(\cC,\cE) = \Sect(\bC,\bfE)$ is left proper combinatorial
for the Reedy model structure.

The rest is done similarly to \cite[Theorem 2.42]{BARLEFT}. As explained there or at \cite[Lemma 2.31]{BALRMSF},
for each object $X \in \bfE(\bc)$ there exists a section $i(X)$ such that
$\PSect(\cC,\cE)(i(X),Y) = \bfE(\bc)(X,Y(\bc))$: indeed, $i(X)(\bc') \cong \coprod_{\alpha: \bc \to \bc'} \alpha_! X$.
For each map $f: \bc \to \bc'$, one has a canonical map $f_! X \to i(X)(\bc')$,
and hence $i(f_! X) \to i(X)$. The required model structure is then obtained by
left Bousfield-localising the Reedy structure on $\PSect(\cC,\cE)$ along the set
$$
\{ i(f_! X) \to i(X) \, \, | \, \, f:\bc \to \bc' \in \sS_\bC, \, X \in G(\bc) \}
$$
with $G(\bc)$ denoting, as in \cite{BARLEFT}, the set of cofibrant homotopy generators of $\bfE(\bc)$. \endproof

Proposition \ref{propdsectasacmc} implies in particular the existence of homotopy
colimits of derived sections, even if the latter are calculated inexplicitly.
To see the same for $\cS$-locally constant derived sections, one can attempt
to use the formalism of right model categories of Barwick \cite{BARRIGHT}.
However, the question of homotopy colimits can be approached differently, using
the language of higher categories.


\section{Higher-categorical aspects}
\subsection{Behaviour with respect to the infinity-localisation}

\label{subsectiondersectcomp}

Our higher-categorical conventions are the same as those adopted in the third section of \cite{BALRMSF},
and we will make use of the results of \cite{BALRMSF} in what follows below.

Let $\cE \to \cC$ be a model opfibration. The subcategory $\cE_f$ of fibrewise
fibrant objects forms a sub-opfibration over $\cC$ with weak equivalence preserving
transition functors. Localising $\cE_f$ along the totality of fibrewise weak equivalences
yields \cite{BALRMSF,HINICH,MGEE1}
a cocartesian fibration in quasicategories $L \cE_f \to \cC$. The goal of this subsection
is to show that $\Sect(\cC,L \cE_f)$ is equivalent to $L \DSect(\cC,\cE)$.

Let us first consider the following diagram:
\begin{equation}
\label{diadersectcomp1}
\begin{diagram}[small]
\cE_f & \lTo & h^*\cE_f & \rTo^S	& \bfE_f \\
\dTo &		&		\dTo	&			&		\dTo \\
\cC & \lTo^h & \bC & \rTo^= &  \bC \\
\end{diagram}
\end{equation}

The functor $S$ was considered in Lemma \ref{simpreplacementrelation}, and is the
restriction to fibrant objects of the right adjoint in the adjunction
$$
S:\cE(c_0) = \Sect([0], c_0^* \cE^\top) \leftrightarrows \Sect([n]^\op, \bc_{[n]}^* \cE^\top)  :ev_0.
$$
Note that the resulting infinity-adjunction (existing by \cite[Corollary 3.31]{BALRMSF} \cite[Theorem 5.1.1]{MGEE1}),
$$
LS:L\cE(c_0)  \leftrightarrows L\bfE(\bc_{[n]}) : Lev_0
$$
has the property that its right adjoint $LS$ is full and faithful, and its essential
image, when restricted $\cE(c_0)_f$,
coincides with the full subcategory $L(\bfE(\bc_{[n]})^{wcart}_f)$ consisting
of fibrant weakly cartesian sections in $\Sect([n]^\op, \bc_{[n]}^* \cE^\top)$, that is, of
those $X: [n]^\op \to \cE^\top$ such that each $X(k)$ is fibrant and each
$X(i) \to \alpha^* X(i-1)$ is a weak equivalence. Indeed, the functors
$$
S:\cE(c_0)_f \leftrightarrows \bfE(\bc_{[n]})^{wcart}_f: ev_0
$$
are inverse weak equivalences \cite{BARKAN} of relative categories.

Given a map $\bc \to \bc'$, the covariant transition functor $\bfE(\bc) \to \bfE(\bc')$, being merely a
pullback, induces $\bfE(\bc)^{wcart}_f \to \bfE(\bc')^{wcart}_f$ and thus an
opfibration $\bfE^{wcart}_f \to \bC$. Combining that with the precedent remarks, using
the functoriality of the infinity-localisation (it can be realised \cite{HINICH} as a functorial fibrant
replacement in marked simplicial sets) we localise the diagram (\ref{diadersectcomp1})
to obtain a diagram of cocartesian fibrations of quasicategories
\begin{equation}
\label{diadersectcomp3}
\begin{diagram}[small]
L\cE_f & \lTo & Lh^*\cE_f & \rTo^{L S}_\sim	& L \bfE^{wcart}_f & \rTo & L\bfE_f \\
\dTo &		&		\dTo	&			&		\dTo	& 	& \dTo \\
\cC & \lTo^h & \bC & \rTo^= &  \bC 	& \rTo & \bC, \\
\end{diagram}
\end{equation}
with the arrow $LS$ being an equivalence by \cite[Proposition 3.1.3.5]{LUHTT}.

The naturality of the Grothendieck construction for quasicategories \cite[Remark 3.1.13]{MGEE1}
together with \cite[Proposition 5.2.3]{MGEE1} imply that the localisation of the opfibration
$\cE_f \to \cC$ is universal, meaning that the map $L h^* \cE_f \to h^* L \cE_f$
is a cocartesian equivalence over $\bC$. On the other hand, since the functor $h: \bC \to \cC$ is an
infinity-localisation along the Segal maps (Proposition \ref{locprop}), we can identify
certain sections over $\bC$ with
the sections over $\cC$. Indeed, given any categorical fibration $\cF \to \cC$, there
is an induced homotopy pull-back diagram
\begin{equation}
\label{diadersectcomp2}
\begin{diagram}[small]
\Fun(\cC,\cF) & \rTo^\sim & \Fun_\sS(\bC,\cF) \\
\dTo 					&						&		\dTo						\\
\Fun(\cC,\cC)	&	\rTo^{\sim}& \Fun_\sS(\bC,\cC)			\\
\end{diagram}
\end{equation}
with $\Fun_\sS(\bC,\cF)$ (and similarly for $\cC$) being the infinity-category of functors $\bC \to \cF$
that send Segal maps of $\bC$ to equivalences of $\cF$. Both vertical maps
in (\ref{diadersectcomp2}) are categorical fibrations, so taking the pull-back
induces an equivalence $\Sect(\cC,\cF) \cong \Sect_\sS(\bC,\cF)$, with the latter
denoting the sub-category of functors $\bC \to \cF$ over $\cC$ sending Segal maps
to equivalences of $\cF$.

Consequently, the infinity-category $\Sect(\cC,L \cE_f)$
is identified with the full subcategory of $\Sect(\bC,h^* L \cE_f)$ consisting of
those $X$ such that for each Segal map $\alpha$ in $\bC$, the map $X(\alpha)$
is cocartesian. For the last statement, one uses the fact that a morphism of $h^* L \cE_f$
is cocartesian iff its image in $L \cE_f$ is such, a direct consequence of
\cite[Remark 2.4.1.4]{LUHTT}.
We conclude the precedent discussion by stating the following.

\begin{lemma}
\label{lmderseccomp1}
Let $\cE \to \cC$ be a model opfibration. Then there is a canonically induced equivalence
between the infinity-category $\Sect(\cC,L \cE_f)$ and the full subcategory of
$\Sect(\bC,L \bfE)$ consisting of those sections $X: \bC \to L \bfE$ such that
\begin{enumerate}
\item For each $\bc \in \bC$, the object $X(\bc)$ belongs to $L \bfE^{wcart}(\bc)$.
\item For each Segal map $\alpha$ in $\bC$, the induced map $X(\alpha)$ is cocartesian.
\end{enumerate}
Furthermore, these two conditions are equivalent to requiring that $X$ sends Segal
maps to \emph{cartesian} maps of $L \bfE$.
\end{lemma}
\proof
Only the last sentence requires proof. For this, note that just as in the 1-categorical
case, the infinity-adjunction $L \bfE(\bc) \rightleftarrows L \bfE(\bc')$ restricts
to an equivalence of quasicategories $L \bfE(\bc)^{wcart} \cong L \bfE(\bc')^{wcart}$ whenever
the map $\bc \to \bc'$ is Segal. Thus for Segal $\alpha:\bc \to \bc'$ and $X$ taking
value in $L \bfE^{wcart}$, the map
$ \bL \alpha_! X(\bc) \to X(\bc')$ is an equivalence iff $ X(\bc) \to \bR \alpha^* X(\bc')$ is such.
\endproof

\begin{prop}
\label{propdersectcomp}
Let $\cE \to \cC$ be a model opfibration. Then the natural comparison functor $L \Sect(\bC, \bfE) \to \Sect(\bC, L \bfE)$
induces an equivalence
$$L \DSect(\cC,\cE) \cong \Sect(\cC,L \cE_f).$$ This equivalence is compatible with the base change.
\end{prop}
\proof
Since $\bfE \to \bC$ is a Quillen presheaf, Theorem 2 (or equally Proposition 3) of \cite{BALRMSF} imply
that $L\Sect(\bC, \bfE) \cong \Sect(\bC, L \bfE)$.
Everything then follows from Lemma \ref{lmderseccomp1} and \cite[Lemma 3.44]{BALRMSF}.
\endproof

A similar result can be proven for the locally constant derived sections.

\begin{thm}
\label{thmcomparisonforlocconstdersect}
Let $\cE \to \cC$ be a model opfibration and $\cS$ a subset of maps of $\cC$.
The infinity-category $\Sect_\cS(\cC,L \cE_f)$ consisting of sections $\cC \to L\cE_f$
that send the maps of $\cS$ to cocartesian maps of $L \cE_f$ is equivalent
to the full subcategory of $\Sect(\bC,\bfE)$ consisting of those $X$ such that
\begin{enumerate}
\item $X$ sends the Segal maps of $\bC$ to cartesian maps of $L \bfE$,
\item $X$ sends the maps of $\bC$ $h_\cC$-projecting to $\cS$, to cocartesian maps of $L \bfE$.
\end{enumerate}
In particular, one has the equivalence $L \DSect_\cS(\cC,\cE) \cong \Sect_\cS(\cC,L \cE_f)$ which
is compatible with the base-change.
\end{thm}
\proof Combine Lemmas \ref{lmderseccomp1}, \ref{lemmalocconstsecconsistency} and \cite[Lemma 3.44]{BALRMSF}. \endproof

As an application of this result, let us see what happens in the presentable setting.
The following is a refinement of \cite[Proposition 5.5.3.17]{LUHTT}.

\begin{lemma}
\label{lmluriepropextended}
Let $p:\cX \to \cY$ be a presentable bicartesian fibration \cite[Definition 5.5.3.2]{LUHTT}
of quasicategories with $\cY$ small,
$\cL, \cR$ two sets of maps such that the covariant transition functors of $p$ along $\cL$
preserve limits.

Let $\Sect_{\cL,\cR}(\cY,\cX)$ denote the full subcategory of $\Sect(\cY,\cX)$
consisting of those sections that send $\cL$ to cocartesian and $\cR$ to cartesian maps in $\cX$.
Then $\Sect_{\cL,\cR}(\cY,\cX)$ is an accessible localisation of the presentable
infinity-category $\Sect(\cY,\cX)$, and hence is presentable itself.
\end{lemma}
\proof
The presentability of $\Sect(\cY,\cX)$ is \cite[Proposition 5.5.3.17]{LUHTT}, and
the proof of the same proposition shows that it suffices to consider $\cY = [1]$ and
prove the accessible localisation condition in the following two cases:

\emph{Case 1:} the inclusion $\Sect_{(\emptyset,\{0 \to 1 \})}([1],\cX) \subset \Sect([1],\cX)$
of cartesian sections into all sections of a presentable fibration $\cX \to [1]$. This
is \cite[Lemma 5.5.3.16]{LUHTT} and the fact that the evaluation at $1 \in [1]$
induces an equivalence $\Sect_{(\emptyset,\{0 \to 1 \})}([1],\cX) \stackrel  \sim \to \cX(1)$.

\emph{Case 2:} the inclusion $\Sect_{(\{0 \to 1 \},\emptyset)}([1],\cX) \subset \Sect([1],\cX)$
of cocartesian sections into all sections of a presentable fibration $\cX \to [1]$ whose
associated covariant transition functor $f_!: \cX(0) \to \cX(1)$ preserves limits.
This can be solved by applying \cite[Proposition 5.4.7.11]{LU} or by the following explicit
argument.
Note that the evaluation at zero induces a categorical equivalence
$\Sect_{(\{0 \to 1 \},\emptyset)}([1],\cX) \to \cX(0)$ (this functor is obtained
by applying cocartesian sections functor to a left marked anodyne map). The association
of a functor $f_!$ to $\cX \to [1]$ provides, by \cite[Definition 5.2.1.1]{LUHTT},
an inverse equivalence $s:\cX(0) \to \Sect_{(\{0 \to 1 \},\emptyset)}([1],\cX)$
such that for $X \in \cX(0)$, $sX(0) = X$ and $sX(1)=f_! X$. Thus given any diagram
$\mathbf Y: I \to \Sect_{(\{0 \to 1 \},\emptyset)}([1],\cX)$ we can assume, that
up to an equivalence, $\mathbf Y= s \mathbf X$, for some $\mathbf X: I \to \cX(0)$.
Taking limits induces the canonical map
$$
f_! \lim_I \mathbf Y(0) = f_! \lim_I \mathbf X \longrightarrow \lim_I f_! \mathbf X = \lim_I \mathbf Y(1);
$$
and by assumption the only non-invertible arrow here is an equivalence. This map
can be viewed as a section of $\cX \to [1]$, and by the dual of (2) of \cite[Proposition 5.1.2.2]{LUHTT}
it serves as a limit of $\mathbf Y$. The category of cocartesian sections of $\cX \to [1]$
is thus closed under limits in $\Sect([1],\cX)$.
By a similar argument, using the fact that $f_!$ is a left adjoint,
the cocartesian section condition is also stable under arbitrary colimits.
Since $\Sect_{(\{0 \to 1 \},\emptyset)}([1],\cX) \to \cX(0)$ is an equivalence, the
category of cocartesian sections is presentable. Thus we have a (co)limit-preserving
inclusion of a presentable category into a presentable category, and such a functor admits
a left adjoint.

The general case follows from $$\Sect_{\cL,\cR}(\cY,\cX) = \cap_{f \in \cL} \Sect_{(\{f \},\emptyset)}(\cY,\cX) \bigcap \cap_{g \in \cR} \Sect_{(\emptyset,\{g \})}(\cY,\cX)$$
and the same reduction argument as in \cite[Proposition 5.5.3.17]{LUHTT}.
\endproof

\begin{corr}
\label{corrlocdsectarepresentable}
Let $\cE \to \cC$ be a model opfibration and $\cS$ a subset of maps of $\cC$. Assume
that the fibres are combinatorial model categories, the transition functors are accessible,
and the transition functors along the maps in $\cS$ preserve limits. Then both
$L \DSect_\cS(\cC,\cE)$ and $\Sect_\cS(\cC,L \cE_f)$ are presentable infinity-categories, with limits and
sufficiently filtered colimits calculated fibrewise.
\end{corr}
\proof Proposition 2.30 of \cite{BALRMSF} implies that the fibres of $\bfE \to \bC$
are combinatorial model categories, and Lemma \ref{lmantisegmapswhenlimitsarepreserved} shows that
the transition functors along the anti-Segal maps that forget the $\cS$-part preserve homotopy limits.
Consequently $L \bfE \to \bC$ is a presentable fibration satisfying the condition of
Lemma \ref{lmluriepropextended} for the sets $\mathscr A \mathscr S_{\cS}$ (of anti-Segal maps of $\bC$ that forget
the $\cS$-labelled part) and $\sS_\bC$ (of Segal maps). It follows from Lemma \ref{lmluriepropextended} and \cite[Lemma 3.44]{BALRMSF} that
the infinity-category
$L \DSect_\cS(\cC,\cE) \cong \Sect_{(\mathscr A \mathscr S_\cS, \sS_\bC)}(\bC, L \bfE)$
is presentable, and hence the same is true for  $\Sect_\cS(\cC,L \cE_f)$. \endproof

\begin{rem}
\label{coalgebrapresentability}
The are other ways to use Theorem \ref{thmcomparisonforlocconstdersect} to conclude
presentability of the infinity-category $L \DSect_\cS (\cC,\cE)$ or of its opposite. Return to the example
of Remark \ref{remexamples} given by the opfibration $\cM^\otimes \to \Fin$, where $\cM = \cN^\op$ for a combinatorial
monoidal model category $\cN$ with all objects cofibrant. Denoting inert maps of $\Fin$ by $In$, it is easy to show
that the infinity-category $\Sect_{In}(\Fin,L \cM_f^\otimes)^\op$ that describes homotopy coalgebra objects in
$\cN$ is presentable. To do so, observe that
$$
\Sect_{In}(\Fin,L \cM_f^\otimes)^\op \cong \Sect_{In}(\Fin^\op,L\cN_c^{\otimes,\top})
$$
which follows from the relation $L \cM_f^{\otimes,\op} \cong L\cN_c^{\otimes,\top}$ (it follows from
\cite[Proposition 5.2.3]{MGEE1} that both these fibrations classify the same infinity-functor) and the
interaction of the functor category with the $\op$-involution.
Then one can further note that due to \cite[Proposition 5.4.7.11 and Remark 5.4.7.14]{LUHTT}
the categories $\Sect_{In}(\Fin^\op,L\cN_c^{\otimes,\top})$ and $\Sect (\Fin^\op,L\cN_c^{\otimes,\top})$
are accessible, and \cite[Proposition 5.1.2.2]{LUHTT} together with the preservation
of colimits along inert maps imply that both $\Sect_{In}(\Fin^\op,L\cN_c^{\otimes,\top})$ and $\Sect (\Fin^\op,L\cN_c^{\otimes,\top})$ are cocomplete. Theorem \ref{thmcomparisonforlocconstdersect}
then implies that the corresponding derived sections $L \DSect_{In}(\Fin,\cM^\otimes)^\op$ that serve
as a rigid model for infinity-coalgebra objects in $\cN$ form a presentable infinity-category.

The same proof works if we replace $\cM^\otimes \to \Fin$ by
any other model opfibration $\cE \to \cC$ whose fibres are opposite-combinatorial,
whose transition functors are opposite-accessible, and that the transition functors along the subset $\cS$
preserve limits. We get that $L \DSect_\cS(\cC,\cE)^\op$ is presentable.
\end{rem}

\subsection{Higher-categorical Segal sections}
\label{highercatsegsections}
The combinatorial trick that leads to derived sections can be reproduced purely in the
higher-categorical setting, and that permits to prove statements like Corollary \ref{corrlocdsectarepresentable}
for more general higher-categorical fibrations.
Let $\cE \to \cC$ be a cocartesian fibration, with $\cE$ an infinity-category and
$\cC$ an ordinary category. In this
subsection we will explicitly construct a cocartesian fibration $\bfE \to \bC$ which
is also cartesian over Segal maps, and identify $\Sect(\cC,\cE)$ with Segal sections inside
$\Sect(\bC,\bfE)$.

\begin{rem}
In the higher-categorical context, we will usually only use the word Segal when
referring to sections, as everything is ``derived'' already.
\end{rem}

The explicit version of the construction that we present uses the
relative nerve functor of Lurie, recalled below. We can also assume that there is a functor $\rE: \cC \to \SSet_+$ taking value in fibrant
objects (that is, $X^\natural$ for $X$ a quasicategory), such that $\cE$ is equivalent
over $\cC$ to the covariant relative nerve of $\rE$.

\begin{lemma}
\label{lemmasegsectinfinitycombinatorics}
Let $\cD: I^\op \to \Cat$ be a diagram of small categories such that each $\cD(i)$
has an initial object $d_i$. Let $\rF_i:\cD(i) \to \SSet_+$ be a family of fibrant-valued
functors, compatible with $\cD$ in the sense that for each $f: i \to i'$, the induced diagram
\begin{diagram}[small,nohug]
\cD(i)  & 		&	\lTo^{f^*} &		& \cD(i')  \\
			&		\rdTo<{\rF_i} & 		&		\ldTo>{\rF_{i'}} &  \\
  & 		&	\SSet_+ &		&  \\
\end{diagram}
commutes up to a canonical isomorphism, with $f^*=\cD(f) $. Then one has a map of projectively fibrant functors in
$\Fun(I,\SSet_+)$ whose value at $i$ is equal to
$$
\rF_i(d_i) \stackrel \sim \longrightarrow \Cart \left( \cD(i)^\op, \int \rF_i \right),
$$
an inverse of the equivalence between the cartesian sections of the \emph{contravariant} relative nerve
$\int \rF_i \to \cD(i)^\op$, and $\rF_i(d_i)$.
\end{lemma}

Note that the assignment $i \mapsto F_i(d_i)$ defines a covariant functor, since
a map $f: i \to i'$ induces a unique map $d_i \to f^* d_{i'}$ in $\cD(i)$. The assignment
$i \mapsto  \Cart \left( \cD(i)^\op, \int \rF_i \right)$ is also covariant, since in view
of the marked version of \cite[Remark 3.2.5.7]{LU}, there is a pullback diagram
\begin{diagram}[size=3em,nohug]
\int F_i 	&		\lTo	& \SWpbk \int  F_{i'} \\
\dTo			&				&		\dTo			\\
\cD(i)^\op & \lTo & \cD(i')^\op	\\
\end{diagram}
that induces the corresponding map of cartesian sections.

\proof
Let us recall the construction of $\int \rF_i \to \cD(i)^\op$. Dualising \cite[Definition 3.2.5.2]{LU},
a map $\Delta^n \to \int \rF_i$ is given by
\begin{enumerate}
\item a map $\sigma :[n] \to \cD(i)^\op$,
\item for each (nonempty) subset $S \subset [n]$ with the minimal element $s$, a map $\Delta^S \to \rF_i(\sigma(s))$,
\item for each such pair $S' \subset S \subset [n]$, the commutativity property of the diagram
\begin{diagram}[small,nohug]
\Delta^{S'} & \rTo & \rF_i(\sigma(s')) \\
\dTo				&			&		\dTo							\\
\Delta^S		&	\rTo	&	\rF_i(\sigma(s)). \\
\end{diagram}
\end{enumerate}
An object of $\int \rF_i$ is thus uniquely determined by a pair consisting of an object
$d \in \cD(i)^\op$ and $x \in \rF_i(d)$. A morphism $(d,x) \to (d',x')$ is the data of a map $\alpha:d' \to d$ in
$\cD(i)$ and a map $x \to \rF_i(\alpha) x'$ in $F_i(d)$, just like one expects from the Grothendieck construction.
The forgetful infinity-functor $\int \rF_i \to \cD(i)^\op$ is a cartesian fibration, with $(d,x) \to (d',x')$ being
cartesian iff $x \to \rF_i(\alpha) x'$ is an equivalence.

To construct the map from the proposition, let us instead construct the adjoint map
$\cD(i)^\op \times \rF_i(d_i) \to \int \rF_i$ over $\cD(i)^\op$. This is done as follows.
A map $\Delta^n \to \cD(i)^\op \times \rF_i(d_i)$ is just a pair of $\sigma: [n] \to \cD(i)^\op$
and $\tau: \Delta^n \to \rF_i(d_i)$. We send it to $ \int \rF_i$ by taking
\begin{enumerate}
\item the same map $\sigma :[n] \to \cD(i)^\op$,
\item for each (nonempty) subset $S \subset [n]$ with the minimal element $s$, take the composition
$$\Delta^S \subset \Delta^n \stackrel \tau \to \rF_i(d_i) \to \rF_i(\sigma(s))$$ to get
a map $\Delta^S \to \rF_i(\sigma(s))$ (here we use that $d_i$ is initial in $\cD(i)$).
\item for each such pair $S' \subset S \subset [n]$, all the squares are commutative in the diagram below,
\begin{equation}
\label{diaplugginginside}
\begin{diagram}[small,nohug]
\Delta^{S'} & \rTo & \Delta^n		&	\rTo^\tau	& \rF_i(d_i)	& \rTo	&		\rF_i(\sigma(s')) \\
\dTo				&			  &		\dTo^=	&						&		\dTo^=		&				&	\dTo							\\
\Delta^S		&	\rTo	&	\Delta^n	&	\rTo_\tau	& \rF_i(d_i)	&	\rTo			&	\rF_i(\sigma(s)), \\
\end{diagram}
\end{equation}
and this implies the commutativity of the outer rectangle.
\end{enumerate}

A map $\Delta^1 \to \cD(i)^\op \times \rF_i(d_i)$ is cartesian iff the image of $\tau: \Delta^1 \to \rF_i(d_i)$
is an equivalence in $\rF_i(d_i)$. It is evidently sent to a cartesian map in $\int \rF_i$. As a result,
the adjoint map $\rF_i(d_i) \to \Sect\left( \cD(i)^\op,\int \rF_i \right)$ factors through the cartesian sections.
Moreover composing with the evaluation map $\Sect\left( \cD(i)^\op,\int \rF_i \right) \to \rF_i(d_i)$ gives
identity. In light of \cite[Corollary 3.3.3.2]{LU} we have thus constructed an equivalence
$$
\rF_i(d_i) \stackrel \sim \longrightarrow \Cart \left( \cD(i)^\op, \int \rF_i \right),
$$
and it remains to check its $I$-naturality.

Given a map $f: i \to i'$, observe that we have the following diagram coming from base-change
\begin{diagram}[size=3em,nohug]
\cD(i)^\op \times \rF_i(d_i) & \rTo & \int \rF_i \\
\uTo 													&			&		\uTo 			\\
\cD(i')^\op \times \rF_i(d_i)	& \rTo & \int \rF_{i'} \\
\dTo													&				&	\dTo>=		\\
\cD(i')^\op \times \rF_i(d_{i'}) & \rTo & \int \rF_{i'} \\
\end{diagram}
with top row projecting to $\cD(i)^\op$ and the bottom square to $\cD(i')^\op$. The middle horizontal
functor is induced by the pullback property of the relative nerve; we can
construct it exactly in such a way that the bottom square commutes. Indeed, for each object $d' \in \cD(i'))$,
we have maps $F_i (d_i) \to F_i (f^* d') = F_{i'}(d')$ that are natural in the $d'$-argument. We can use them to construct
$\cD(i')^\op \times \rF_i(d_i)	\to \int \rF_{i'}$ similarly to we did above at (\ref{diaplugginginside}) and verify everything by hand.
Taking the cartesian sections then gives the diagram
\begin{diagram}[size=3em,nohug]
\rF_i(d_i)			&	\rTo		&	\Cart(\cD(i)^\op,\cD(i)^\op \times \rF_i(d_i)) & \rTo & \Cart\left(\cD(i)^\op,\int \rF_i \right) \\
\dTo<=					&					&	\uTo 													&			&		\uTo 			\\
\rF_{i}(d_i)		&	\rTo		&	\Cart(\cD(i')^\op,\cD(i')^\op \times \rF_i(d_i))	& \rTo & \Cart \left(\cD(i')^\op,\int \rF_{i'} \right) \\
\dTo						&					&	\dTo													&				&	\dTo>=		\\
\rF_{i'}(d_{i'})&	\rTo		&	\Cart(\cD(i')^\op,\cD(i')^\op \times \rF_{i'}(d_{i'})) & \rTo & \Cart \left(\cD(i')^\op,\int \rF_{i'} \right). \\
\end{diagram}
One has a canonical isomorphism of infinity-categories $$\Cart(\cD(i)^\op,\cD(i)^\op \times \rF_i(d_i)) \cong \Fun_{lc}(\cD(i)^\op,\rF_i(d_i)),$$ with the right hand side denoting functors taking all arrows to equivalences. Under this identification the maps $\rF_i(d_i) \to 	\Cart(\cD(i)^\op,\cD(i)^\op \times \rF_i(d_i))$
correspond to constant diagram functors, and this proves the commutativity of the diagram above. Taking
its outer rectangle proves the sought-after naturality.
\endproof

The next lemma will provide adjoints along Segal maps.

\begin{lemma}
\label{lemmainfinitysegalmaps}
Let $\cF \to [n]$ be a cartesian fibration.
Let $\alpha: [m] \subset [n]$ be an interval inclusion such that $\alpha(m)=n$ and
for each $i \in [n]$ not in the image of $\alpha$, the fibre $\cF_i$ has an initial
object. Then the pull-back functor $\alpha^*: \Sect([n],\cF) \to \Sect([m],\cF)$
has both right and left adjoints $\alpha_*$ and $\alpha_!$. The right adjoint furthermore factors
through the full subcategory of $\Sect([n],\cF)$ consisting of those sections that
are cartesian over the complement of the image of $\alpha$.
\end{lemma}

\proof Begin first with the case of $i: \{ 1 \} \subset [1]$. Given a cartesian fibration $\cF \to [1]$
and associating to it a functor as in \cite[Definition 5.2.1.1]{LU} provides, among other things, us with a functor
 $i_*: \cF_1 \to \Cart([1],\cF) \subset \Sect([1],\cF)$ such that $i^* i_* = id$. The diagram of \cite[Lemma 5.5.3.15]{LU}
 (or the relative Kan extension formalism) then tells us that $i_*$ is right adjoint to $i^*$.

For the left adjoint, note that the initial object $\emptyset \in \cF_0$ is also initial
in the whole $\cF$ thanks to \cite[Lemma 5.2.3.4]{LU}. Thus for any $X \in \cF_1$
there exists a map $\emptyset \to X$ in $\cF$ that projects to $0 \to 1$; we can view this map
as a section. Using \cite[Corollary 4.3.2.16 and Proposition 4.3.2.17]{LU} we see that
there thus exists a functor $i_!: \cF_1 \to \Sect([1],\cF)$ left adjoint to $i^*$
such that $i^* i_! = id$.

For the general case, it will be enough to show the existence of adjoints with required properties
for the inclusions $\beta:[n-1] \subset [n]$. This inclusion induces the pushout-product
map $\bar \beta: [1] \cup_{[0] } [n-1] \to [n]$ where $[0]$ includes as $0$ in $[n-1]$ and as $1$ in $[1]$.
Using \cite[Lemma 2.1.2.3]{LU} for $\emptyset \to [0]$ and the left anodyne map $[0] \subset [n-1]$,
we conclude that the map $\bar \beta$ is inner anodyne and thus $$\Sect([n],\cF) \to \Sect([1],\cF) \times_{\cF_1} \Sect([n-1],\cF)$$
is an equivalence of infinity-categories. It moreover respects the partial cartesian property:
a section $[n] \to \cF$ is cartesian over $0 \to 1$ iff its image in the fibred product has a
cartesian $\Sect([1],\cF)$-component.
It will thus suffice to construct the adjoints for the functor
$$ \Sect([1],\cF) \times_{\cF_1} \Sect([n-1],\cF) \to  \cF_1 \times_{\cF_1} \Sect([n-1],\cF) \cong \Sect([n-1],\cF),$$
but we can do this simply by taking fibred products of the adjoint triple $i_! \dashv i^* \dashv i_*$
constructed before with
$\Sect([n-1],\cF) \to \cF_1$ and remembering that mappings spaces of fibred products are fibred products of mapping
spaces. \endproof

\begin{rem}
Implicit in the proof of Lemma \ref{lemmainfinitysegalmaps} is the fact that given a
cartesian (or cocartesian) fibration $\cE \to \cB$, the functor $\Sect(-,\cE)$
is left Quillen with respect to the Joyal model model structure over $\cB$, as
follows from, for example, \cite[Remark 3.1.4.5 and Proposition 3.1.5.3]{LU}.
As a consequence all the fibred products appearing in the proof are also homotopy fibred
products.
\end{rem}

We are now ready to formulate our higher-categorical variation of derived sections.
Note that given a cocartesian fibration $\cE \to \cC$, we will denote by $\cE^\top \to \cC^\op$
its transpose or dual, as per usual. It can be obtained through the Grothendieck construction
procedure, or with the help of \cite{BARGLANAR}. The following proposition summarises the
Segal section construction:

\begin{prop}
Let $\cC$ be a small category and $\cE \to \cC$ a cocartesian fibration. Denoting
$h_\cC:\bC \to \cC$ the initial element functor as before, there exists
a sequence
$$
h^*_\cC \cE \stackrel \cong \longrightarrow \bfE^{cart} \longrightarrow \bfE
$$
of cocartesian fibrations and functors over $\bC$, with the first functor being a
cocartesian equivalence and the second also preserving cocartesian arrows. The fibre
of this diagram at $\bc:[n] \to \cC$ is
$$
\cE(c_0) \stackrel \cong \longrightarrow \Cart([n]^\op,\bc^* \cE^\top) \longrightarrow \Sect([n]^\op,\bc^* \cE^\top)
$$
with the first functor being an inverse of the evaluation equivalence. The transition
functors between the diagrams for different fibres are induced by the transition functors of $
\cE$ and the evident restrictions of (cartesian) sections.

Finally, the cocartesian fibrations $\bfE \to \bC$ and $\bfE^{cart} \to \bC$ are also cartesian along Segal maps.
The composition $h^*_\cC \cE \to \bfE$ induces an equivalence between
$\Sect(\cC,\cE)$ and the subcategory $\Sect_{\sS}(\bC,\bfE)$ consisting of those
sections $X:\bC \to \bfE$ that send Segal maps to cartesian maps. This equivalence is compatible
with the base change along 1-functors.
\end{prop}

\proof
We can replace $\cE \to \cC$ by $\int_\cC \rE \to \cC$, the covariant relative nerve
of a (projectively fibrant) functor $\rE: \cC \to \SSet_+$. The transpose fibration
$\cE^\top \to \cC^\op$ can be then replaced by $\int_{\cC^\op} \rE \to \cC^\op$,
the contravariant relative nerve of $\rE$. Given $\bc_{[n]}: [n] \to \cC$, an object of $\bC$,
let us denote $\rE_\bc := \rE \circ \bc_{[n]}$. We can then apply
Lemma \ref{lemmasegsectinfinitycombinatorics} for the diagram $\bC^\op \to \Cat, \bc_{[n]} \mapsto [n]$
and the system of functors $\rE_{\bc_{[n]}}$.
The result is a diagram of functors from $\bC$ to $\SSet_+$ whose value at $\bc_{[n]}:[n] \to \cC$
is
$$
\rE(c_0) \stackrel \cong \longrightarrow \Cart \left([n]^\op, \int \rE_{\bc_{[n]}} \right)
\longrightarrow \Sect \left([n]^\op, \int \rE_{\bc_{[n]}} \right)
$$
with the contravariant relative nerve underneath. All what remains is to take the covariant relative nerve
and apply Lemma \ref{lemmainfinitysegalmaps} for the extra adjoints.

The section identification then proceeds as usual. Given a bicartesian fibration $\cM \to [1]$
representing an equivalence $\cM_0 \cong \cM_1$, a section $X: [1] \to \cM$ is cartesian iff
it is cocartesian. Thus any section $X: \bC \to \bfE^{cart}$ sends Segal maps to cartesian maps
iff it sends them to cocartesian maps.
Since $h_\cC: \bC \to \cC$ is
an infinity-localisation along the Segal maps, we have a diagram
$$
\Sect(\cC,\cE) \stackrel \cong \longleftarrow \Sect_\sS (\bC,\bfE^{cart}) \longrightarrow \Sect_\sS(\bC,\bfE)
$$
with first arrow an equivalence, and the $\sS$-index denoting the Segal section condition.
However, the second arrow is also an equivalence, as it follows from the description of right
adjoints along $(c_0 \to ... \to c_n) \to c_0$, as given in Lemma \ref{lemmainfinitysegalmaps},
that any Segal section $\bC \to \bfE$ factors through $\bfE^{cart}$.
\endproof

\begin{corr}
\label{corrlocconsitinfintyderivedsect}
Let $\cI$ be a subset of morphisms of $\cC$. Then the infinity-category $\Sect_\cI (\cC,\cE)$
of $\cI$-cocartesian sections is naturally identified with the infinity-category
$\Sect_{\cI-loc,\sS}(\bC,\bfE)$ consisting of Segal sections that send the $\cI$-coloured
anti-Segal maps, $$(c_0 \stackrel \cI \to ... \stackrel \cI \to c_k \to ... \to c_n) \longrightarrow (c_k \to ... \to c_n)$$
(see Definition \ref{deflocconstdersect}), to cocartesian maps.
\end{corr}
\proof
A map of $\bC$ that projects to $\cI$ can be factored as a composition of an $\cI$-coloured
anti-Segal map and of an initial element preserving map. Any Segal section is automatically
cocartesian along the latter type of maps.
\endproof

Let us now use the higher-categorical Segal sections to see the presentability
of section categories.

\begin{thm}
\label{thmpresentabilitygeneral}
Let $\cE \to \cB$ be a cocartesian fibration over a small infinity category $\cB$, such
that each fibre $\cE(b)$ is presentable and each transition functor is accessible.
Let $\cS \subset \cB_1$ be a subset of maps such that each transition functor $f_!$ along
$f \in \cS$ preserves limits (thus satisfies the right adjoint functor theorem).

Then the infinity-category of $\cS$-cocartesian sections $\Sect_\cS(\cB,\cE)$ is presentable, with limits
and (sufficiently large) filtered colimits calculated fibrewise.
\end{thm}
This result is not found in \cite[5.5]{LU}, and can be used to abstractly conclude the
existence of colimits and adjoints of algebras over infinity-operads and operator categories.

\proof We already know that $\Sect_\cS(\cB,\cE)$ is accessible, with limits and filtered colimits of sections
calculated fibrewise,
as follows from \cite[Proposition 5.1.2.2 and Corollary 5.4.7.17]{LU}.
Choose an infinity-localisation $F:\cC \to \cB$ surjective on objects and morphisms,
where $\cC$ is an ordinary category. This can be done by taking for example the simplex
category of $\cB$, see \cite{STEVENSON}. Denote by $F^*\cS$ the subset of maps
of $\cC$ that $F$-project to $\cS$ or to $\cB$-equivalences.
The induced functor $\Sect_\cS(\cB,\cE) \to \Sect_{F^* \cS}(\cC,F^* \cE)$
is then seen to be an equivalence of infinity-categories: if $X \cong F^* Y$ and $X$
is $F^* \cS$-cocartesian, then $Y$ is $\cS$-cocartesian since any map in $\cS$ has a preimage in $F^* \cS$.

Having a cocartesian fibration $\cE \to \cC$ over an ordinary category and a subset
of maps $\cI$, we can identify $\Sect_\cI(\cC,\cE)$ with $\Sect_{\cI-loc,\sS}(\bC,\bfE)$
as per Corollary \ref{corrlocconsitinfintyderivedsect}. Furthermore, the fibres
of $\bfE \to \bC$ are sections of pullbacks of $\cE^\top \to \cC^\op$. The latter is a
cartesian fibration in presentable categories and accessible transition functors.
By \cite[Proposition 5.1.2.2 and Corollary 5.4.7.17]{LU} in this case we have that
the fibres $\bfE(\bc)$ are presentable, with colimits calculated objectwise, which
in turn implies that the transition functors of $\bfE \to \bC$ preserve colimits.
This in turn implies that $\bfE \to \bC$ is a presentable fibration. All that remains is
to check that the $\cI$-coloured anti-Segal maps lead to limit-preserving transition functors.

The following will suffice. Let $\cF \to [n]$ be a cartesian fibration in presentable infinity-categories
and accessible functors. Let $\omega: [m] \subset [n]$ be
an interval inclusion such that $\omega(0) = 0$. Assume that for each $i\geq m$, the transition
functors $\cF_{i+1} \to \cF_i$ preserve limits. Then the pullback
functor $\omega^*: \Sect([n],\cF) \to \Sect([m], \cF)$ also preserves limits.

To prove this statement, fix an inclusion $[n - m] \subset [m]$ that is complement to $\omega$.
With the help of \cite[Lemma 2.1.2.3]{LU}, it will be enough to show that the functor
$\Sect([m],\cF) \times_{\cF_m} \Sect([n-m], \cF) \to \Sect([m],\cF)$ preserves limits.
Given a fibration in complete infinity-categories and limit-preserving functors,
the category of sections is also complete, with limits calculated fibrewise: this is
\cite[Proposition 5.4.7.11]{LU} together with \cite[Remark 5.4.7.13]{LU}, but can also be proven by hand.
Thus the category $\Sect([n-m], \cF)$ has fibrewise limits. The category $\Sect([m],\cF)$
is presentable, and the functor of evaluation at $m$ preserves all limits and colimits by
Lemma \ref{lemmainfinitysegalmaps}. Thus the (homotopy) fibred product $\Sect([m],\cF) \times_{\cF_m} \Sect([n-m], \cF)$
is the limit of a diagram of complete categories and limit-preserving functors.
By \cite[Proposition 5.4.7.11 and Remark 5.4.7.13]{LU} we once again have that this fibre product
is complete with limits calculated in each category separately, and thus preserved by any projection.

Returning to the proof of the theorem, since the transition functors along $\cI$
preserve limits, we can apply the precedent discussion and then conclude everything using
Lemma \ref{lmluriepropextended}.
\endproof

\section{Resolutions}


Given a functor $F: \cD \to \cC$ and a model opfibration $\cE \to \cC$, we get an induced pullback functor $\bF^*: \PSect(\cC,\cE) = \Sect(\bC,\bfE) \to \Sect(\bD,\bfE) = \PSect(\cD,\cE)$ on the categories of presections. This functor trivially preserves weak equivalences, and also the condition of being a derived section. We thus get the functor
$$
\bF^*: \Ho \DSect(\cC,\cE) \to \Ho \DSect (\cD, \cE)
$$
on the level of localisations.
In some cases, the functor $\bF^*$ is full and faithful, and its essential image is easy to characterise.

\begin{opr}
	\label{weakfibresdefinition}
	For a functor $F: \cD \to \cC$ and $\bc_{[n]} = c_0 \to ... \to c_n$ of $\bC$, denote by $\cD(\bc_{[n]})$ the category
	\begin{itemize}
		\item with objects being pairs of $\bd_{[n]} = d_0 \to ... \to d_n$ and of a commutative diagram
		\begin{diagram}[small,nohug]
		Fd_0 & \rTo & ... & \rTo &  Fd_n \\
		\dTo^\cong & 		& ... 		&		& \dTo^\cong \\
		c_0 & \rTo & ... & \rTo &  c_n \\
		\end{diagram}
		so that the vertical maps are isomorphisms,

		\item with morphisms given by commutative diagrams
		\begin{diagram}[small]
		d_0 & \rTo & ... & \rTo &  d_n \\
		\dTo & 		& ... 		&		& \dTo \\
		d'_0 & \rTo & ... & \rTo &  d'_n \\
		\end{diagram}
		such that for each $0 \leq i \leq n$, the diagram
		\begin{diagram}[small,nohug]
		F d_i & \rTo & F d'_i \\
		\dTo<\cong & 	& \dTo^\cong \\
		c_i & \rTo^= & c_i \\
		\end{diagram}
		commutes.
	\end{itemize}
\end{opr}

The categories $\cD(c_0 \to ... \to c_n)$ are extensions of the notion of an essential fibre of a functor.

\begin{opr}
	\label{resolutiondefinition}
	$\,$
	\begin{enumerate}
		\item A functor $F: \cD \to \cC$ is a \emph{resolution} if for each $\bc_{[n]} \in \bC$, the category $\cD(\bc_{[n]})$ is contractible (that is, has a contractible nerve).

		\item A functor $F: \cD \to \cC$ is a \emph{right resolution} if
		\begin{itemize}
			\item for each $c \in \bC$ over $[0] \in \Delta$, the category $\cD(c)$ is contractible, and
			\item for each $f: c' \to c$ in $\bC$ over $[1] \in \Delta$ and $d \in \cD(c)$, the subcategory $F(f,d) \subset \cD(c' \stackrel f \to c)$ given by the (strict) fibre of $\cD(c' \stackrel f \to c) \to \cD(c)$ over $d$, is contractible.
		\end{itemize}
		\item A functor $F: \cD \to \cC$ is a \emph{left resolution} if
		\begin{itemize}
			\item for each $c \in \bC$ over $[0] \in \Delta$, the category $\cD(c)$ is contractible, and
			\item for each $f: c' \to c$ in $\bC$ over $[1] \in \Delta$ and $d \in \cD(c)$, the subcategory $F(d',f) \subset \cD(c' \stackrel f \to c)$ given by the (strict) fibre of $\cD(c' \stackrel f \to c) \to \cD(c')$ over $d'$, is contractible.
		\end{itemize}
	\end{enumerate}
\end{opr}

\begin{lemma}
	If $F:\cD \to \cC$ is a right or left resolution, then $F$ is also a resolution.
\end{lemma}
\proof
We prove the right part, the left part being dual, Inductively, assume we have proven the resolution property for each $\bc'_{[k]}$ with $0 \leq k < n$. Then for an object $\bc_{[n]}= c_0 \stackrel f \to c_1 \to ... \to c_n$ we have the associated functor  $\cD( c_0 \stackrel f \to c_1 \to ... \to c_n) \to  \cD(c_1 \to ... \to c_n)$. This is an opfibration over a contractible category, with fibres equivalent to $F(f,d)$ for some $d \in \cD(c_1)$. Quillen's Theorem A implies then the contractibility of $\cD(\bc_{[n]})$. \endproof

\begin{lemma}
	\label{whenprefibisaresolution}
	\label{whenpreopfibisaresolution}
	If $F: \cD \to \cC$ is a prefibration (and, by convention, an isofibration) with contractible fibres, then it is a right resolution.

	Dually, if $F: \cD \to \cC$ is a preopfibration (and, by convention, an isofibration) with contractible fibres, then it is a left resolution.
\end{lemma}
\proof
Since $F$ is an isofibration, the categories $\cD(c)$ and $\cD(c' \to c)$ are equivalent to their strict analogues: the strict fibre $F^{-1}(c)$ and the category of arrows $d' \to d$ with $F(d' \to d)$ equal to $c' \to c$. It is then easy to see that the fibres of $ \cD(c' \to c) \to \cD(c)$ have terminal objects for a prefibration, and the fibres of $\cD(c' \to c) \to \cD(c')$ have initial objects for a preopfibration, and hence are contractible. Quillen's Theorem A, again, implies the result. \endproof

This example motivates the intuition behind left (respectively right) resolutions as covariant (respectively contravariant) families of contractible homotopy types, represented in categories, indexed by the base category $\cC$.

\begin{lemma}
	\label{adjunctiongivesresolutions}
	If $p: \cD \rightleftarrows \cC: i$ is an adjunction and $i$ is full and faithful, then $p$ is a resolution. An equivalence of categories is a resolution.
\end{lemma}
\proof
Every fibre $\cD(c_0 \to ... \to c_n)$ has a terminal object given by $i c_0 \to ... \to i c_n$ (note that $pi (c_0) \to ... \to pi (c_n)$ is isomorphic to $c_0 \to ... \to c_n$).
\endproof


In \cite[Key Lemma]{HINICH}, it is shown that a resolution $F: \cD \to \cC$ exhibits
the category $\cC$ as a higher-categorical localisation of $\cD$ with respect to
$F$-isomorphisms. Let us reproduce the argument here.

\begin{prop}[Hinich's Lemma]
Let $F: \cD \to \cC$ be a resolution. Then $F$ induces an equivalence $L \cD \cong \cC$,
where we infinity-localise $\cD$ along the set $F^* Iso(\cC)$ of all maps of $\cD$ that are sent
by $F$ to isomorphisms in $\cC$.
\end{prop}
\proof
It is enough to show that
$$
(\cD,F^* Iso(\cC)) \longrightarrow (\cC,Iso(\cC))
$$
is a weak equivalence of relative categories. Following \cite{BARKAN},
one has to show that for each $[n]$, the induced functor
$$
\Fun_{F^* Iso(\cC)}([n],\cD) \to \Fun_{Iso(\cC)}([n],\cC)
$$
is a homotopy equivalence; here the sub-index of $\Fun$ means that we take
only those natural transformations that belong to ($F$-)isomorphisms.
The comma-fibres of the functor above are given by
$\cD(c_0 \to ... \to c_n)$. \endproof

The comparison between ordinary sections and derived sections carried out in
\ref{subsectiondersectcomp} and the Key Lemma of Hinich readily implies
that derived sections on $\cC$ are identified with the derived sections on $\cD$
which are locally constant along the $F$-isomorphisms. It turns out however that proving
the same statement internally to the language of derived sections, without passing
by the comparison of \cite[Theorem 2]{BALRMSF}, leads to some interesting simplicial combinatorics.

Denote as usual by $\DSect_\cS(\cC,\cE)$ the category of $\cS$-locally constant derived sections.
We prove the following result independently of the comparison carried out in \ref{subsectiondersectcomp}:
\begin{prop}
\label{resolutionthm}
	Let $F: \cD \to \cC$ be a resolution, $\cS \subset \cC$ a subset, and $\cE \to \cC$ a model opfibration. Then the functor  $\bF^*: \DSect_\cS(\cC,\cE) \to \DSect_{F^* \cS} (\cD, \cE)$ is a weak equivalence of relative categories, meaning the
	equivalence of induced simplicial (or higher-categorical) localisations.
\end{prop}


\begin{rem}
For many uses it is already sufficient to know the weaker statement, that the functor
$\bF^*: \Ho \DSect_\cS(\cC,\cE) \to \Ho \DSect_{F^* \cS} (\cD, \cE)$ is an equivalence.
However proving such a weaker statement often gives more for free.
To verify that a given relative functor is a weak equivalence, one can often use the following
criterion. Let $F: \cM \rightleftarrows \cN:G$ be a Quillen adjunction with $G$ preserving all weak equivalences,
and assume that there is
a full subcategory $\cM_0 \subset \cM$ closed under the weak equivalences of $\cM$. If
$G$ factors through $\cM_0$ and the homotopy category adjunction $\bL F: \Ho \cM \rightleftarrows \Ho \cN: \bR G$ restricts to an equivalence
between $\Ho \cM_0$ and $\Ho \cN$, then the same is true for the induced higher-categorical adjunction
existing thanks to \cite{CISBOOK,MGEE1}, and so $G: \cM \to \cM_0$ gives an equivalence after
applying the higher-categorical localisation functor. Thus in practice one is left working
with Quillen adjunctions and their restrictions, as we shall indeed do later in this section.
\end{rem}

\subsection{Relative comma objects}
\begin{ntn}
For a $\Delta$-indexed category $\pi:\cX \to \Delta^\op$, the notation $x \twoheadrightarrow y$ means that the map is a degeneracy, that is the underlying map $\pi(y) \twoheadrightarrow \pi(x)$ is a surjection in $\Delta$.
\end{ntn}

\begin{opr}
	Let $F:\cX \rightarrow \cZ \leftarrow \cY:G$ be a diagram of $\Delta$-indexed categories. The associated \emph{$\Delta$-relative comma object} $F \sslash G$ is the full subcategory of the ordinary comma-category $F / G$ consisting of all triples $(x,y, F(x) \twoheadrightarrow G(y))$ where the map $F(x) \twoheadrightarrow G(y)$ is a degeneracy.
\end{opr}

By definition, $F \sslash G$ comes with projections to both $\cX$ and $\cY$. We will also use the notation $F \bbslash G$ to denote $G \sslash F$. If one of the functors is an identity, we will write, as is customary for comma categories, $F \sslash \cY$ instead of $F \sslash id_\cY$ and the like. To know more about the relative comma objects, the following may be of use.

\begin{lemma}
Let $[1]$ be the usual arrow category. Then the full subcategory $Arr_s(\Delta)$ of the arrow category $\Fun([1],\Delta)$ consisting of surjective arrows, is naturally equipped with a Reedy category structure.
The natural source and target projections from $Arr_s(\Delta)$ to $\Delta$ are compatible
with the Reedy structure.
\end{lemma}
\proof Any map between two surjective arrows in $\Delta$ is represented by a commutative square. Using the injection-surjection factorisation system, we get the diagram
\begin{diagram}[small]
[n] & \rOnto &  [n''] & \rInto & [n'] \\
\dOnto &	&	\dTo &		& \dOnto \\
[m] & \rOnto &  [m''] & \rInto & [m'] \\
\end{diagram}
in which it is easy to check that the middle vertical arrow is also a surjection.
To see the Reedy category structure, set $deg([n] \twoheadrightarrow [m]) = n+m$. \endproof

\begin{corr}
\label{reedystructureoncomma}
Let $F:\cX \rightarrow \cZ \leftarrow \cY:G$ be a diagram of $\Delta$-indexed categories. Then
the category $F \bbslash G$ is a $Arr_s(\Delta)$-indexed category, hence a Reedy category,
and both functors $F \bbslash G \to \cX$ and $F \bbslash G \to \cY$ are compatible with the
Reedy structure.
\end{corr}
\proof Immediate from Proposition \ref{indexedcatinheritance}. \endproof

\begin{lemma}
For a $\Delta$-indexed functor $F: \cY \to \cX$, the projection $F \bbslash \cX \to \cX$ is a Reedy functor which is moreover an opfibration over face maps of $\cX$.
\end{lemma}

\proof Corollary \ref{reedystructureoncomma} implies that the projection is Reedy. Now, fix a face map $f: x \to x'$ and $o \in F \bbslash \cX$ with $p(o)=x$. Projecting these data to $\Delta$, we find ourselves with an injection $\phi: [n] \hookrightarrow [m]$ and a surjection $\omega: [k] \twoheadrightarrow [m]$ representing the object $o$. Form the pullback square
\begin{diagram}[small]
[l]\SEpbk			&	\rInto			&	[k]		\\
\dOnto			&			&	\dOnto>\omega \\
[n]		&	\rInto^\phi			&	[m]		\\
\end{diagram}
which exists for this particular configuration of arrows, and use the fact that $F \bbslash \cX$ is discretely fibred over $Arr_s(\Delta)$ to uniquely reconstruct the map $o \to f_! o$ in $F \bbslash \cX$.
\endproof

In the proof above, note that if $f$ is a degree-lowering map, then so is $o \to f_! o$.

\begin{opr}
A functor $G: \cR \to \cR'$ between Reedy categories is called \emph{right-compatible},
if it preserves the Reedy structure and in addition for each $x \in \cR$,
the induced map $Mat(x) \to Mat(Gx)$ is initial relative to $\cR'$, by which we mean that for any functor $X$
from $\cR'$ to a complete category, there is a (naturally induced) isomorphism
$$
\lim_{Mat(Gx)} X|_{Mat(Gx)} \cong \lim_{Mat(x)} G^*X|_{Mat(x)}.
$$
\end{opr}

\begin{corr}
\label{pxisrightcompat}
The functor $p_\cX:  \cO:= F \bbslash \cX \to \cX$ is a right-compatible functor of Reedy categories.
\end{corr}

\proof The opfibration property over face maps yields us an adjunction between the matching categories
$$
L:Mat(p_\cX o) \rightleftarrows Mat(o):R
$$
where $R$ is the canonical functor induced by $p_\cX$. The functor $L$, being a (fully faithful) left adjoint, is initial. The composition of $L$ with the natural projection $p:Mat(o) \to \cO \to \cX$ coincides with the projection $q:Mat(p_\cX o) \to \cX$. We thus can observe that for any functor $X: \cX \to \cM$ valued in a complete category, there is the following sequence of isomorphisms
$$
\lim_{Mat(o)} p^* X \cong \lim_{Mat(p_\cX o)} L^* p^* X \cong \lim_{Mat(p_\cX o)} q^* X
$$
which implies right-compatibility.  \endproof

The main result of this subsection is the interpretation of the higher categorical localisation
property of resolutions in the context of simplicial replacements:

\begin{prop}
\label{propertiesofcommaresolution}
	Let $F: \cD \to \cC$ be a functor. Then for the functor $\bF: \bD \to \bC$, the projection $p_\bC: \bF \bbslash \bC \to \bC$ is a right-compatible Reedy functor. Moreover, if $F$ is an isofibration and a resolution, then the functor $p_\bC$ is a left resolution.
\end{prop}

A proof is given below. We can thus produce the following technical definition.

\begin{opr}
	Let $F: \cY \to \cX$ be a $\Delta$-indexed functor. $F$ is called a \emph{resolution} if the associated functor $F \bbslash \cX \to \cX$ is a left resolution.
\end{opr}
\begin{rem}
Due to $\cX$ being a Reedy category any functor $\cY \to \cX$ is automatically an isofibration.
\end{rem}
In the course of the proof, we shall need an auxiliary lemma about cofinal maps.

\begin{lemma}
\label{lemmaintegrationcofinality}
Consider the diagram
\begin{diagram}[small,nohug]
\cE & \rTo^f & \cF \\
\dTo<p & 			& \dTo>q \\
\cD & \rTo_g & \cC \\
\end{diagram}
and assume that $p,q$ are opfibrations, and both $g$ and all $f_x: \cE(x) \to \cF(g(x))$
(induced by taking fibres over $x \in \cD$) are homotopy cofinal. Then the map $f$
is homotopy cofinal.
\end{lemma}
\proof
For a diagram $X: \cF \to \Top$, we have that the map $\hocolim_\cE f^* X \to  \hocolim_\cF X$
decomposes in $\Ho \Top$ as follows:
$$
\hocolim_\cE f^* X \cong \hocolim_{\cD} \bL p_! f^* X \cong \hocolim_\cD g^* \bL q_! X
$$
$$
\cong \hocolim_\cC \bL q_! X
\cong \hocolim_\cF X.
$$
Here the first and the last isomorphisms are due to the properties of homotopy left Kan extensions
as homotopy adjoints, the third one uses homotopy cofinality of $g$ and the second one
follows from the fact that the base change $\bL p_! f^* \to g^* \bL q_!$ being an
isomorphism, which follows from the requirement of the homotopy cofinality condition on the fibres
and the fact that one calculates $(\bL p_! Y)(d)$ as $\hocolim_{\cE(d)} Y|_{\cE(d)}$,
and similarly for $\bL q_!$. \endproof

\proof[Proof of Proposition \ref{propertiesofcommaresolution}.]
An object of $\bF \bbslash \bC$ is represented by a surjection $\sigma:[k_0 + ... + k_n] \twoheadrightarrow [n]$, an object
$$
\bd_{[k_0 + ... + k_n]} = d_0^0 \to ... \to d^{k_0}_0 \to d^{0}_1 \to ... \to d^{k_n}_n
$$
of maps in $\cD$ such that $p_\bC (\bd_{[k_0 + ... + k_n]})$ is equal to
$$
\sigma^* \bc_{[n]} = c_0 \to ... \to c_0 \to c_1 ... \to c_n
$$
with each $c_i$ appearing $k_i +1$ times in a row. Sending $\bd_{[k_0 + ... + k_n]}$ to $d^{0}_0 \to d^{0}_1 \to ... \to d^{0}_n$ is seen to produce a functor
$$
\tau_{\bc_{[n]}}:(\bF \bbslash \bC)(\bc_{[n]}) \to  \cD(\bc_{[n]}).
$$ Effectively, we are taking the beginning of each sub-division of $\bd_{[k_0 + ... + k_n]}$,
which looks similar to (but is different from) the projection from the simplicial
replacement of a category to the category itself.

We now prove that the functors $\tau_{\bc_{[n]}}$ are homotopy cofinal, hence induce
homotopy equivalences. For the case of a single object, $\bc = c_0$, we see that we are simply comparing the category $\cD(c)$ and its simplicial replacement. And it is well known that
for any category $\cX$, the initial element functor $\mathbb X \to \cX$ is homotopy
cofinal.

By induction, consider the following diagram
\begin{diagram}[small]
\left(\bF \bbslash \bC \right)(\bc_{[n]}) & \rTo^{\tau_n} &  \cD(\bc_{[n]}) \\
\dTo		&		& \dTo \\
\left(\bF \bbslash \bC \right)(\bc_{[n-1]}) & \rTo^{\tau_{n-1}} &  \cD(\bc_{[n-1]}) \\
\end{diagram}
with $\bc_{[n-1]} = c_1 \to ... \to c_{n}$ and both vertical functors given by natural projections. The bottom arrow $\tau_{n-1}$ is homotopy cofinal, both vertical arrows are opfibrations, and the restriction of $\tau_n$ on the fibres of the left arrow gives, again, the standard functor between a category and its simplicial replacement. Thus $\tau_n$ is homotopy cofinal by Lemma \ref{lemmaintegrationcofinality}.

To continue, we also need to study the fibres of the projection

\begin{equation}
\label{heavyproofproj}
 (\bF \bbslash \bC)(\bc'_{[k]} \to \bc_{[n]}) \to  (\bF \bbslash \bC)(\bc'_{[k]}).
\end{equation}
Fix a morphism $s:\bc'_{[k]} \to \bc_{[n]}$ in $\bC$ and an object
$(\bd,\bc'_{[k]},f)$ in $(\bF \bbslash \bC)(\bc'_{[k]})$. Denote by
$Fibre(s,\bd,\bc'_{[k]},f)$ the fibre of the projection (\ref{heavyproofproj})
over $(\bd,\bc'_{[k]},f)$.

An object of $Fibre(s,\bd,\bc'_{[k]},f)$ is, by definition, an object
$(\bd',\bc_{[n]},f')$ and a morphism $s':\bd \to \bd'$ such that
$s \circ f = f' \circ \bF(s')$. Since $\bC$ is discretely opfibred over $\Delta$,
the problem is completely defined by its image in $\Delta$. The category
$Fibre(s,\bd,\bc'_{[k]},f)$ does not, thus, depend on the exact detail of the
categories $\cD,\cC$, so we can replace them with one-object categories.
Effectively, we are given a map $g:[n] \to [k]$ and a surjective map $h:[m] \twoheadrightarrow [k]$, and we consider triples $[m'],h',g'$, with $h':[m'] \twoheadrightarrow [n]$ surjective, $g':[m'] \to [m]$ arbitrary, and $h \circ g' = h' \circ g$.

Factor $[n] \stackrel g \to [k]$ as a surjection and an injection, $[n] \stackrel{g_s}{\to} [n''] \stackrel{g_i}{\to} [k]$, and observe that we can take pullbacks of surjections along $g_i$, with results being surjections. We thus see that we are studying the category of possible diagrams
\begin{diagram}[small, nohug]
[m'] & \rOnto & [m'' \SEpbk ]  & \rInto & [m] \\
\dOnto					&		& \dOnto &				& 	\dOnto>h \\
[n] & \rOnto^{g_s} &	[n''] & \rInto^{g_i} & [k] \\
\end{diagram}
where the whole right (pullback) square and $g_s$ are fixed. The data of $[m'] \twoheadrightarrow [n]$ is equivalent to an object of $\Delta^{n+1}$. Specifying a compatible map $[m'] \to [m'']$ then gives us a functor $L: (\Delta^{n+1})^\op \to \Set$: there are no non-trivial morphisms between two different liftings. Moreover, if we denote by $g_{s,*}: \Delta^{n+1} \to \Delta^{n''+1}$ the post-composition functor, we see that $L \cong g_{s,*}^* S$, where $S: (\Delta^{n'' +1})^\op \to \Set$ is the functor represented by $[m''] \twoheadrightarrow [n'']$.

It will thus suffice to prove the following. Consider any surjection $g:[n] \twoheadrightarrow [k]$, and the induced functor $g_*:\Delta^{n+1} \to \Delta^{k+1}$. Then we need to show that for any representable functor $S:(\Delta^{k+1})^\op \to \Set$, the $n+1$-fold simplicial set $(g_*)^*S$ is contractible. By induction on $k$, it suffices to consider the case $k=0$. Take the diagonal embedding $\delta:\Delta \to \Delta^{n+1}$. Then $|\delta^*(g_*)^*S|$ is equivalent to $|(g_*)^*S|$, so it suffices to prove that for any $X: \Delta^\op \to \Set$, one has a homotopy equivalence $|i_{n+1}^*X| \cong |X|$, where $i_{n+1} = g_* \circ \delta:\Delta \to \Delta$. Explicitly, $i_{n+1}$ acts exactly as $n+1$-fold edgewise subdivision functor, and $|i_{n+1}^*X|$ is homotopically equivalent (actually homeomorphic) to $|X|$ for any simplicial set $X$.
\endproof

\subsection{Sections over relative comma objects}

In the following, it will be useful to axiomatise the relevant properties of the functor $\bfE \to \bC$.

\begin{opr}
\label{segalbifibdef}
	Let $\cX$ be a $\Delta$-indexed category. A \emph{model Segal bifibration} over $\cX$ is a bifibration $\bfE \to \cX$ with the following properties.
	\begin{enumerate}
		\item The bifibration $\bfE \to \cX$ is a Quillen presheaf: the fibres $\bfE(x)$ are model categories, and for each $f: x \to y$, the adjunction $f_!: \bfE(x) \rightleftarrows \bfE(y): f^*$ is a Quillen pair.
		\item The following base change condition holds. For a commutative square
		\begin{diagram}[small]
		x & \rTo^f & y \\
		\dTo<\alpha & &  \dTo>\beta \\
		z & \rTo_g & t \\
		\end{diagram}
with vertical arrows Segal, and horizontal arrows projecting to initial element preserving maps in $\Delta$, the induced natural transformation $\bL f_! \bR \alpha^* \to \bR \beta^* \bL g_!$ is an isomorphism.
	\end{enumerate}
\end{opr}

By \cite[Theorem 1]{BALRMSF} we have that the category of sections $\Sect(\cX,\bfE)$ carries a Reedy model structure.

\begin{opr}
Let $\bfE \to \cX$ be a model Segal bifibration, as described above. A section $X: \cX \to \bfE$ is called \emph{Segal} if for any Segal map $f:x \to y$ the induced morphism
$$
X(x) \to \bR f^* X(y)
$$
is an isomorphism in $\Ho \bfE(x)$.
\end{opr}

We denote by $\Sect_{\sS}(\cX,\bfE)$ the corresponding full subcategory.

\begin{lemma}
\label{stabilityalongdegeneracies}
Let $X$ be a Segal section of $\bfE \to \cX$. Then for each degeneracy $s:x \twoheadrightarrow y$ of $\cX$, both induced maps $X(x) \to \bR s^* X(y)$ and $\bL s_! X(x) \to X(y)$ are isomorphisms in $\Ho \bfE(x)$ and $ \Ho \bfE(y)$ respectively.
\end{lemma}

\proof Similar to Lemma \ref{derivedsectionsdontseedegenerations}. \endproof

For a functor $F: \cY \to \cX$, we would like to study what happens when we lift Segal sections from $\cX$ to $F \bbslash \cX$ via the natural projection $p_\cX$. In general, the  latter is an opfibration over face maps, which implies the following.

\begin{prop}
\label{pushforwardformulaleftres}
Let $\bfE \to \cX$ be a model bifibration over a $\Delta$-indexed category $\cX$, and $p: \cO \to \cX$ be a right-compatible Reedy functor. Then there is an induced Quillen adjunction
$$
p_!: \Sect(\cO,p^*\bfE) \rightleftarrows \Sect(\cX,\bfE):p^*
$$
where the right adjoint is the pullback functor.

Furthermore, if $p$ is a left resolution, then the pushforward $\bL p_!$ admits a simplified expression
$$
\bL p_! X(x) \cong \bL \colim_{\cO(x)} X|_{\cO(x)},
$$
suitably functorial in $x$.
\end{prop}

\proof By definition, the functor $p$ identifies the matching objects computed for $o$ in $\cO$ with those of its image $p(o)$. This permits to verify that $p^*$ preserves Reedy fibrations and trivial Reedy fibrations. Its left adjoint, $p_!$, is not hard to compute. Given $x \in \cX$ and a section $Y: \cO \to \bfE$, denote by $p/x$ the usual comma category fibre consisting of pairs $(o,p(o) \to x)$; we then have
$$
p_! Y(x) \cong \colim_{p/x} Res_x Y|_{p / x}
$$
with $Res_x: Sect(p/x, \bfE) \to \bfE(x)^{p/x}$ being the restriction functor induced by the opfibration structure on $\bfE$. The category $p/x$ is Reedy with the same latching objects as in $\cO$, and the functor $Res_x$, as can be checked, preserves (trivial) Reedy cofibrations. This implies the same expression for the derived functor,
$$
\mathbb L p_! Y(x) \cong \mathbb L \colim_{p/x} \mathbb L Res_x Y|_{p / x}.
$$
It remains to see that the inclusion functor $i:\cO(x) \to p/x$ is homotopy cofinal in the case of a left resolution. For this, we consider the categories $z \backslash i$, where $z = (o,p(o) \to x)$ is an object of the comma category $p/x$. One can see that the category $z \backslash i$ is the category of commutative squares
\begin{diagram}[small]
o			&	\rTo			&	o'		\\
\dMapsto			&			&	\dMapsto \\
p(o)			&	\rTo			&	x		\\
\end{diagram}
where the upper left object and the bottom arrow are fixed. It is the same category as the fibre over $o$ of the projection $\cO(p(o) \to x) \to \cO(p(o))$, and it is contractible. Which in turn implies that we can restrict the homotopy colimits from $p/x$ to $\cO(x)$. \endproof

\begin{corr}
\label{corollaryleftres}
In the situation of Proposition \ref{pushforwardformulaleftres}, if the functor $p:\cO \to \cX$ is a left resolution, then the pullback functor $p^*$ factors through the full subcategory $\Sect_{loc}(\cO,p^*\bfE)$ consisting of those sections $Y:\cO \to p^*\bfE$ which send fibrewise maps in $\cO$ to weak equivalences in the fibres of $\bfE$. Moreover
$p^*: \Sect(\cX,\bfE) \to \Sect_{loc}(\cO,p^*\bfE)$ is a weak equivalence of relative categories.
%
%
\end{corr}

\proof
For $Y = p^* X$, we see that
$$
\bL p_! Y(x) \cong \bL \colim_{\cO(x)} X(x) \cong X(x)
$$
since the category $\cO(x)$ is contractible. Conversely, given $Y$ which is homotopically constant on the fibres, we have that for any $o \in \cO(x)$, which exists since $\cO(x)$ is nonempty, the natural map $Y(o) \to \bL \colim_{\cO(x)} Y$ is an isomorphism in $\Ho \bfE(x)$ \cite{CIS}. However, $\bL \colim_{\cO(x)} Y \cong (p_! Y)(x) \cong (p^* p_! Y)(o)$. \endproof


Returning to the case of the relative comma category $F \bbslash \cX$ for $F: \cY \to \cX$, we note that, besides $p_\cX$, there is another projection functor $p_\cY:F \bbslash \cX \to \cY$, which is also fairly special.

\begin{lemma}
\label{propertiesofpy}
The functor $p_\cY$ is a discrete opfibration admitting a fully faithful right adjoint $i_\cY: \cY \to F \bbslash \cX$ such that $F= p_\cX i_\cY$. The functor $p_\cY$ is also a right-compatible Reedy functor.
\end{lemma}
\proof The right adjoint is easy to see, and is given by sending $y \in \cY$ to $i_\cY y = (y, F(y) \stackrel = \twoheadleftarrow F(y))$, its fully-faithfulness is then obvious. Its opfibration property is also apparent, as the diagram
\begin{diagram}[small]
F(y)			&	\rTo			&	F(y')		\\
\uOnto			&			&	  \\
x		&	 		&			\\
\end{diagram}
can be completed to
\begin{equation}
\label{auxdiag1}
\begin{diagram}[small]
F(y)			&	\rTo			&	F(y')		\\
\uOnto			&			& \uOnto	  \\
x		&	 \rTo		&	x'		\\
\end{diagram}
\end{equation}
using the face-degeneracy factorisation system inherited from $\Delta$. As in Corollary \ref{pxisrightcompat}, the lift (\ref{auxdiag1}) allows us to conclude that $p_\cY$ is right-compatible. \endproof

\begin{prop}
\label{liftingtocommafromdomain}
Let $\bfE \to \cX$ be a model Segal bifibration over a $\Delta$-indexed category $\cX$, and $F: \cY \to \cX$ be a $\Delta$-indexed functor. Then the functor $p_\cY: F \bbslash \cX \to \cY$ induces a Quillen adjunction
$$
i_\cY^*=p_{\cY,!}: \Sect( F \bbslash \cX, p^*_\cY F^*\bfE) \rightleftarrows \Sect(\cY,F^* \bfE):p_\cY^*
$$
where $i_\cY$ is the functor from Lemma \ref{propertiesofpy}. The functor $p_\cY^*$ is a weak equivalence between the relative category  $\Sect(\cY,F^* \bfE)$ and the full subcategory of $\Sect( F \bbslash \cX, p^*_\cY F^*\bfE)$ consisting of those $Y$ such that for each $p_\cY$-fibrewise map $f$ in $F \bbslash \cX$, the induced map $Y(f)$ is a weak equivalence.
\end{prop}
\proof
The adjunction $i^*_\cY \dashv p_\cY^*$ is easy to verify, and the right-compatibility of $p_\cY$ makes its pullback right Quillen. We now study the essential image of $p_\cY$. Note that for each $y$, the fibre $y \bbslash \cX$ of the opfibration $p_\cY$ contains a final object $i_\cY y$, and is in particular contractible. In a way similar to Corollary \ref{corollaryleftres}, we can now see that the essential image of $\mathbb R p_\cY^*$ consists of sections homotopically constant on fibres.
\endproof

\subsection{Comparing the Segal sections}


\begin{opr}
\label{segalmapsincomma}
A map $s:o \to o'$ of $F \bbslash \cX$ is \emph{Segal} if both its projections $p_\cX s$ and $p_\cY s$ are Segal maps in $\cX$ and $\cY$ respectively.
\end{opr}


\begin{lemma}
\label{segalmapslemma1}
Let
\begin{equation}
\label{segalmapslemma}
\begin{diagram}[small]
[l]			&	\lTo^\alpha			&	[k]		\\
\dOnto			&			&	\dOnto \\
[n]		&	\lTo^\beta			&	[m]		\\
\end{diagram}
\end{equation}
be a diagram in $\Delta$.
\begin{enumerate}
\item If the map $\alpha$ is a left (respectively, right) interval inclusion, then we can factor (\ref{segalmapslemma}) as
\begin{diagram}[small]
[l]			&	\lTo^\alpha			&	[k] & \lTo^= & [k]\\
\dOnto			&			&	\dOnto & & \dOnto \\
[n]		&	\lTo^\gamma			&	[m'] & \lTo^\zeta & [m]		\\
\end{diagram}
with $\gamma \zeta = \beta$, $\gamma$ being a left (respectively, right) interval inclusion, and $\zeta$ a surjection.
\item If the map $\beta$ is a left (respectively, right) interval inclusion, then we can factor (\ref{segalmapslemma}) as
\begin{diagram}[small]
[l]			&	\lTo^\epsilon			&	\SWpbk [k'] & \lTo^\mu & [k]\\
\dOnto			&			&	\dOnto & & \dOnto \\
[n]		&	\lTo^\beta			&	[m] & \lTo^= & [m]		\\
\end{diagram}
with $\alpha = \epsilon \mu$, and $\epsilon$ being a left (respectively, right) interval inclusion, and the left square a pull-back.
\end{enumerate}
\end{lemma}
\proof The second factorisation is evident. For the first, consider all the elements contained in the image $[k] \stackrel \alpha \longrightarrow [l] \longrightarrow [n]$, they will form an ordered subset included as a left (respectively, right) interval. The map $\zeta$ is then forced to be a surjection. \endproof

\begin{lemma}
\label{segalmapslemma2}
Given a surjection $[n] \twoheadrightarrow [m]$ and a left (right) interval inclusion $\alpha: [k] \to [l]$, there exists a completion to a square
\begin{diagram}[small]
[l]			&	\lTo^\alpha			&	[k]		\\
\dOnto			&			&	\dOnto \\
[n]		&	\lTo^\beta			&	[m]		\\
\end{diagram}
with horizontal arrows being left (right) interval inclusions. Same with $\alpha$ replaced by $\beta: [m] \to [n]$
\end{lemma}
\proof Clear. \endproof

For a functor $F: \cY \to \cX$ and a model Segal bifibration $\bfE \to \cX$, there are two bifibrations over the relative comma category which one can obtain using pullback operations: $p_\cX^* \bfE \to F \bbslash \cX$ and $p_\cY^* F^* \bfE \to F \bbslash \cX$. The bifibration structure induces a Quillen adjunction
$$
S:\Sect(F \bbslash \cX, p^*_\cX \bfE) \rightleftarrows \Sect(F \bbslash \cX, p_\cY^* F^* \bfE):T.
$$
Over an object $(y, F(y) \stackrel f \twoheadleftarrow x)$, the functor $S$ amounts to post-composing with $f_!: \bfE(x) \to \bfE(F(y))$, and the functor $T$ -- with its right adjoint $f^*$.

\begin{prop}
\label{bifibchangeprop}
The Quillen adjunction $S \dashv T$ induces an adjoint equivalence
$$
\bL S:\Ho \Sect_{Seg}(F \bbslash \cX, p^*_\cX \bfE) \rightleftarrows \Ho \Sect_{Seg}(F \bbslash \cX, p_\cY^* F^* \bfE) : \bR T$$
where $\Ho \Sect_{Seg} (F \bbslash \cX, p^*_\cX \bfE)$ denotes the full subcategory of sections $Y:F \bbslash \cX \to  p^*_\cX \bfE$ for which the map $Y(o) \to \bR t^* Y(o')$ induced from a Segal map $t: o \to o'$ is an isomorphism (similarly for $ p_\cY^* F^* \bfE$).
\end{prop}
\proof The restriction to Segal sections is possible due to base change property of the bifibration in question (Definition \ref{segalbifibdef}), as the adjoints are essentially defined
using the transition functors.
We shall analyse both the unit and the counit of the adjunction restricted to Segal sections.

Take a Segal section $Y:F \bbslash \cX \to  p^*_\cX \bfE$. Over $o=(y, F(y) \stackrel f \twoheadleftarrow x)$, the map $Y(o) \to \bR T \bL S Y (o)$ is isomorphic to $Y(o) \to \bR f^* \bL f_! Y(o)$. Denote by $x_0$ and $y_0$ the first vertices of $x$ and $y$; the evident map $o=(y, F(y) \stackrel f \twoheadleftarrow x) \to (y_0, F(y_0)=x_0)=o_0$ is then Segal. If we denote
\begin{diagram}[small]
x & \rOnto^f & F(y) \\
\dTo<{\pi_x} &     & \dTo>{\pi_y} \\
x_0 & \rOnto^{f_0} & F(y_0) \\
\end{diagram}
then, using the Segal condition, we see that the map $$Y(o) \to \bR f^* \bL f_! Y(o)$$ is isomorphic to $$\bR \pi_x^* Y(o_0) \to \bR f^* \bL f_! \bR \pi_x^* Y(o_0).$$ Finally, using the base change from Definition \ref{segalbifibdef}, we see that the latter map is isomorphic to $$\bR \pi_x^* Y(o_0) \to \bR f_0^*\bR \pi_y^* Y(o_0),$$ which is merely a differently written identity map.

The counit of the adjunction is treated similarly.
\endproof

\begin{prop}[Resolutions for $\Delta$-indexed categories]
	\label{resolutiontheoremrough}
Let $\bfE \to \cX$ be a model Segal bifibration over a $\Delta$-indexed category $\cX$, and $F: \cY \to \cX$ be a $\Delta$-indexed functor. If $F$ is a resolution, then
$$
F^* : \Sect_{Seg}(\cX,\bfE) \to \Sect_{Seg}(\cY,F^* \bfE)
$$
factors through those Segal sections $X: \cY \to F^* \bfE$ such that $\bR p_\cY^* X$ is homotopically constant on the fibres of $p_\cX: F \bbslash \cX \to \cX$, and is a weak equivalence onto the subcategory of such $X$.
\end{prop}

\proof
Corollary \ref{corollaryleftres} and Proposition \ref{liftingtocommafromdomain} can be enhanced using Lemmas \ref{segalmapslemma1} and \ref{segalmapslemma2}. For the projection $p_\cX: F \bbslash \cX \to \cX$, we get the following result: the homotopical essential image of
$$
p_\cX^*: \Sect_{Seg}(\cX,\bfE) \to \Sect(F \bbslash \cX, p^*_\cX \bfE)
$$
consists of those sections $Y$ for which
\begin{enumerate}
\item the map $Y(o) \to \bR t^* Y(o')$ induced from a Segal map $t: o \to o'$ of $F \bbslash \cX$ is an isomorphism,
\item the usual condition applies, that is, homotopical constancy on the fibres of $p_\cX$.
\end{enumerate}
The necessary observation is that $\bL p_{\cX,!} Y$ is Segal, as one can verify by lifting Segal maps from $\cX$ to $F \bbslash \cX$ using Lemma \ref{segalmapslemma2}.

For the projection $p_\cY$, we get, formally, that
the homotopical essential image of
$$
p_\cY^*: \Sect_{Seg}(\cY,F^*\bfE) \to \Sect(F \bbslash \cX, p^*_\cY F^* \bfE)
$$
consists of those sections $Z$ for which
\begin{enumerate}
\item the map $Z(o) \to \bR t^* Z(o')$ induced from a Segal map $t: o \to o'$ of $F \bbslash \cX$ is an isomorphism,
\item the usual condition applies, that is, homotopical constancy on the fibres of $p_\cY$.
\end{enumerate}

Now we need to pass between two different bifibrations. Proposition \ref{bifibchangeprop} implies that no information is lost by applying $T$ and $S$ to the sections of interest.

Let $Y: F \bbslash \cX \to p_\cX^* \bfE$ be a Segal section constant on the fibres of $p_\cX$. We then see that $\bL S (Y)$ is a Segal section which is automatically constant on the fibres of $p_\cY$: these fibres are $y \bbslash \cX$, with only possible maps being induced from degeneracies in $\cX$, and any derived section coming from $\cX$ is constant along degeneracies (Lemma \ref{stabilityalongdegeneracies}). Thus $\bL S(Y)$ is contained in the essential image of $p^*_\cY$.

For the converse, we see that in order for $Z:F \bbslash \cX \to p_\cY^* F^* \bfE$ to belong to the essential image of $p^*_\cX$, we need exactly that $\bR T (Z)$ is homotopically constant on the fibres of $p_\cX$. Both precedent observations and the fact that $\bR T \bR p_\cY^* X$ is homotopically constant on the fibres of $p_\cX$ iff $\bR p_\cY^* X$ is such imply the theorem. \endproof

The first corollary is a result which is, perhaps surprisingly, not obvious.

\begin{corr}
\label{corequivalencesareequivalences}
Let $F: \cD \to \cC$ be an equivalence of categories and $\cE \to \cC$ a model opfibration. Then the induced pullback
functor $\bF^* : \DSect(\cC,\cE) \to \DSect(\cD,\cE)$ is an equivalence of relative categories.
\end{corr}
\proof Denote by $\cI$ the category whose objects are triples $(d,c, F(d) \cong c)$.
Both natural projections $\pi_\cC: \cI \to \cC$ and $\pi_\cD: \cI \to \cD$ are equivalences (hence resolutions)
and isofibrations; $\pi_\cD$ admits a canonical inverse equivalence $i_\cD$ such that $\pi_\cC \circ i_\cD = F$.
The functor $\pi_\cC$ also admits a section $i_\cC$ that corresponds to choosing an equivalence
inverse to $F$. Denote by $\pi_\bC,\pi_\bD,i_\bC, i_\bD$ the corresponding simplicial replacement functors.

A map $f$ in $\cI$ is an isomorphism iff it is true for either $\pi_\cC f$ or $\pi_\cD f$.
This together with Proposition \ref{resolutiontheoremrough} implies that the
pullbacks $\pi_\bC^* : \DSect(\cC,\cE) \to \DSect(\cI,\pi_\cC^* \cE)$ and
$\pi_\bD^* : \DSect(\cD,F^*\cE) \to \DSect(\cI,\pi_\cD^* F^* \cE)$ are weak equivalences
of relative categories. The same is thus true for $i_\bC^*$ and $i_\bD^*$.

Denote by $[[n]]$ the category that has two copies of $[n]$ inside of it,
denoted $0 \to ... \to n$ and $0' \to ... \to n'$, and exactly one isomorphism
$i \cong i'$ for each $0 \leq i \leq n$, with the obvious commutativity identities.
Any object $(c,d,Fd \cong c)$ of $\cI$ yields a map $[[0]] \to \cC$; similarly,
it is true that given an object $\sigma$ of the simplicial replacement $\bI$,
we have an associated functor $\sigma : [[n]] \to \cC$ whose first row looks like
$F d_i$ and the second row like $c_i$. We can use this notation to define
$\tilde \bfE (\sigma) = \Sect([[n]], \sigma^* \cE^\top)$ and form the covariant Grothendieck
construction to get $\tilde \bfE \to \bI$. Evaluating sections at the top or the bottom row
induces functors $\pi_\bD^* \bF^* \bfE \stackrel q \leftarrow \tilde \bfE \stackrel p \rightarrow \pi_\bC^* \bfE$ over $\bI$
that are furthermore seen to be equivalences. We can promote $\tilde \bfE \to \bI$
to a model Segal bifibration, such that the equivalences $q$ and $p$ preserve all
parts of the (fibrewise) model structure. It makes sense to speak of Segal
sections of $\tilde \bfE \to \bI$.

Observe now that we have the following diagram:
\begin{diagram}[small,nohug]
													&								&	\Sect_\sS(\bI, \tilde \bfE)			&									&												\\
													& \ldTo^{p_*} 			&			&	\rdTo^{q_*} 			&												\\
\DSect(\cI,\pi_\cC^* \cE) & 							&			&				&		\DSect(\cI,\pi_\cD^* F^* \cE)		\\
\dTo<{i_\bC^*}  & 							&			&				&	 \dTo>{i_\bD^*} 	\\
\DSect(\cC, 				 \cE) & 							&	\rTo^{\bF^*}		&				&		\DSect(\cD,  F^* \cE)		\\
\end{diagram}
This diagram commutes up to an isomorphism: a section $X: \bI \to \tilde \bfE$ can
be evaluated on those $\sigma: [[n]] \to \cC$ that look like
\begin{diagram}[small, nohug]
F d_0 & \rTo & ... & \rTo & F d_n \\
\dTo<=	&		&				&	&  \dTo>= \\
F d_0 & \rTo & ... & \rTo & F d_n \\
\end{diagram}
creating an isomorphism between the restriction of $X(\sigma)$ to the top and the bottom
rows. Moreover $p_*,q_*,i_\bC^*,i_\bD^*$ are weak equivalences; it thus follows that
$\bF^*$ is one as well.  \endproof

\proof[Proof of Proposition \ref{resolutionthm}.] Given an arbitrary resolution $F: \cD \to \cC$,
we can factor it as $F = \tilde F G$, where $\tilde F$ is an isofibration and a resolution,
and $G$ is an equivalence. Applying Proposition \ref{resolutiontheoremrough} for $\tilde \bF$
and Corollary \ref{corequivalencesareequivalences} for $G$, we get that
$\DSect(\cC,\cE) \to \DSect_{F^* Iso_\cC}(\cD,F^* \cE)$ is a weak equivalence of relative
categories, with $F^* Iso_\cC$ denoting the $F$-isomorphisms. Given any subset $\cS$ of maps of $\cC$,
if we ever have a map
$$
(c_0 \stackrel \cS \to ... \stackrel \cS \to c_k \to c_{k+1} \to ... \to c_n) \longrightarrow (c_k \to c_{k+1} \to ... \to c_n)
$$
in $\bC$, then the resolution property of $F$ allows us to find a map in $\bD$ that looks like
$$
(d_0 \stackrel {F^*\cS}{ \to} ... \stackrel {F^*\cS}{ \to} d_k \to d_{k+1} \to ... \to d_n) \longrightarrow (d_k \to d_{k+1} \to ... \to d_n)
$$
and such that its $\bF$-image covers the $\bC$-map in question. This observation allows us to see
that any $F^* \cS$-constant derived section $X$ is the $\bF^*$-image of a $\cS$-constant derived section.
\endproof

\end{document}